\documentclass[12pt]{amsart}
\usepackage{amssymb}
\usepackage{amsbsy}

\usepackage{euscript}
\let\cal\EuScript
\let\script\EuScript
\let\epsilon\varepsilon

\def\BMO/{\ifmmode{\mathop{\mathrm{BMO}}}\else {\normalfont BMO}\fi}
\def\BMOr/{\ifmmode{\mathop{\mathrm{BMOr}}}\else {\normalfont BMOr}\fi}
\def\BMOAr/{\ifmmode{\mathop{\mathrm{BMOAr}}}\else {\normalfont BMOAr}\fi}
\def\BMOA/{\ifmmode{\mathop{\mathrm{BMOA}}}\else {\normalfont BMOA}\fi}
\def\bmo/{\ifmmode{\mathop{\mathrm{bmo}}}\else {\normalfont bmo}\fi}

\makeatletter
\def\eqalign#1{\null\,\vcenter{\openup\jot\m@th
  \ialign{\strut\hfil$\displaystyle{##}$&$\displaystyle{{}##}$\hfil
      \crcr#1\crcr}}\,}
\def\eqalignbot#1{\null\,\vbox{\openup\jot\m@th
  \ialign{\strut\hfil$\displaystyle{##}$&$\displaystyle{{}##}$\hfil
      \crcr#1\crcr}}\,}
\def\eqaligntop#1{\null\,\vtop{\openup\jot\m@th
  \ialign{\strut\hfil$\displaystyle{##}$&$\displaystyle{{}##}$\hfil
      \crcr#1\crcr}}\,}

\let\over\@@over

\@mparswitchfalse

\@addtoreset{equation}{section}

\makeatother

\def\eqref#1{\textnormal{(\ref{eq#1})}}

\newtheorem{lemma}{Lemma}[section]
\newtheorem{proposition}[lemma]{Proposition}
\newtheorem{theorem}[lemma]{Theorem}
\newtheorem{corollary}[lemma]{Corollary}

\newtheorem{unnumberedtheorem}{Theorem}

\newtheorem{unnumberedcorollary}[unnumberedtheorem]{Corollary}

\theoremstyle{definition}
\newtheorem{remark}{Remark}
\newtheorem{definition}{Definition}

\newbox\provedbox
\renewenvironment{proof}
  {\setbox\provedbox\hbox{$\square$}\trivlist\item[]{\bf Proof.}}
  {\ifvoid\provedbox\else\hproved\fi\endtrivlist}
\def\proved{\ifmmode\eqno{\box\provedbox}\else\hproved\fi}
\def\hproved{\unskip\nobreak\hfil\penalty50\hskip.5em\hbox{}\nobreak\hfil
  \box\provedbox{\parfillskip=0pt\finalhyphendemerits=0\par}}

\def\tri{\|\mskip-1.8mu|}           
\let\ell=l              
\let\<=\langle
\let\>=\rangle
\def\dist{\mathop{\mathrm{dist}}\nolimits}
\def\codim{\mathop{\mathrm{codim}}}
\def\Re{\mathop{\mathrm{Re}}}

\def\N{\mathbb N}
\def\C{\mathbb C}
\def\D{\mathbb D}
\def\T{\mathbb T}
\def\funnyT{\cal I}
\def\bee{b}

\begin{document}

\title[Nehari--AAK Theory for Hankel Operators on the Torus]
{Two Distinguished Subspaces of Product \BMO/\\
and the Nehari--AAK Theory\\ for Hankel Operators on the Torus}
\author{Mischa Cotlar}
\address{Mischa  Cotlar\\
Facultad de Ciencias\\
Universidad Central de Venezula\\
Caracas 1040, Venezuela}
\email{cs@@scs.howard.edu}

\author{Cora Sadosky}
\address{Cora Sadosky\\
Department of Mathematics\\
Howard University\\
Washington, DC 20059, USA}
\email{mcotlar@@dino.conicit.ve}

\thanks{Sadosky was partially supported by NSF grants DMS-9205926, INT-9204043
and GER-9550373, and her visit to MSRI is supported by NSF grant
DMS-9022140 to MSRI}

\subjclass{47B35, 42B20}

\begin{abstract}
In this paper we show that
the theory of Hankel operators in the torus $\T^d$, for $d > 1$,
presents striking differences with that on the circle $\T$, starting with
bounded Hankel operators with no bounded symbols. Such differences are
circumvented here by replacing the space of symbols $L^\infty (\T)$ by
$\BMOr/(\T^d)$, a subspace of product \BMO/, and the singular numbers of Hankel
operators by so-called sigma numbers.
This leads to versions of the Nehari--AAK
and Kronecker theorems, and provides conditions
for the existence of solutions of product Pick problems
through finite Pick-type matrices.  We give geometric
and duality characterizations
of \BMOr/, and of a subspace of it,
\bmo/, closely linked with $A_2$ weights.  This
completes some aspects of the theory of \BMO/ in product spaces.
\end{abstract}

\maketitle

\section*{Introduction}

This paper deals with the extension of the classical theory of bounded
Hankel operators in the circle $\T$ to (big) Hankel operators in
the torus $\T^d$, for $d > 1$.  Some crucial results in the one-variable
theory, involving the notions of $L^\infty$ symbols and the singular
numbers of the operators, cannot have, as stated, meaningful
extensions to the torus.  This difficulty can be overcome by
introducing so-called \BMOr/ symbols and sigma numbers of Hankel
operators.  To explain what changes are to be made in dimension $d>1$,
we recall some basic features of the theory in $\T$.

Each function $\phi \in L^2(\T^d)$ gives rise to a Hankel
operator $\Gamma _\phi $, and $\phi $ is called a symbol for the operator.  In
the case $d = 1$,
these operators are closely related to the space \BMO/, since,
by the Nehari theorem [N], a Hankel operator $\Gamma $ is bounded if and only if
$\Gamma 1\in \BMO/$, and if and only if $\Gamma  = \Gamma _\varphi $ with
$\varphi \in L^\infty $, while $\phi \in  \BMO/$ implies
$\Gamma _\phi $ bounded with $\|\Gamma _\phi \| = \|\phi \|_{\BMO/}$.  In turn,
the Helson--Szeg\H{o} theorem [HS] relates \BMO/ to the boundedness of the
Hilbert transform in $L^2(\mu )$, for $\mu $ a given measure on the circle $\T$.

The Nehari theorem gives the distance of a bounded function
$\varphi $ to the space $H^\infty (\T)$ as the norm of the Hankel operator
$\Gamma _\varphi $, and the theorem of Adamjan, Arov and Krein (AAK) refines
this by giving its distance to $H^\infty (\T) + R_n$ (where $R_n$ is the space
of rational functions with $n$ poles in the disk) as the singular number $s_n$
of the operator, or, equivalently, as the distance of the operator to those
Hankel operators of finite rank $n$ [AAK].  

From the Beurling characterization of the invariant subspaces of
$H^2(\T)$ of finite codimension, it follows that a Hankel operator
$\Gamma$ is of finite rank $n$ if and only if $\Gamma=\Gamma_\phi$
with $\phi=\bar bh$, where $h\in H^\infty$ and $b$ is a Blaschke
product with $n$ zeros at $z_1,\dots,z_n$. If this is the case,
the operator $\Gamma_\phi$ is closely related to a model operator in a
finite subspace of $H^2$, so that its norm $\|\Gamma_\phi\|$ equals
that of a finite $n\times n$ matrix explicitly given in terms of the
$z_k\!$'s and $\phi(z_k)$'s: the Pick matrix.  One of the main
applications of the Nehari theorem is that it provides a condition for
the existence of solutions of the Pick interpolation problems in terms
of the norm of an associated Hankel operator $\Gamma_\phi$ of finite
rank, thus yielding the classical Pick condition in terms of Pick
matrices.

The basic properties of $\BMO/(\T)$ can be deduced in a unified way [ACS]
through a generalized Bochner theorem, which includes also the results of
Nehari and Helson--Szeg\H{o}.  The
extension of this theorem to several dimensions led in [CS2] and [CS3] to an
extension of the Nehari theorem to $\T^d$, for $d > 1$, in terms of a class of
symbols that we called \BMOr/ (for ``restricted'' \BMO/).  The extension of the
Helson--Szeg\H{o} theorem to several dimensions was given in [CS1], in terms
of a subspace of product $\BMO/ = \BMO/(\T^d)$ (defined in [ChF1]), that here we
call \bmo/ (for ``small'' \BMO/).

Section 1 gives
some basic properties of these subspaces of product \BMO/,
starting with the continuous proper inclusions
$$
L^\infty (\T^d) \subset  \bmo/ \subset  \BMOr/ \subset  \BMO/(\T^d).
$$
The preduals of \bmo/ and \BMOr/ are determined, providing counterparts of the
duality result of Chang and Fefferman in product domains [ChF2].

As a corollary of the duality result for \BMOr/, in Section 2 it is shown that,
when $d > 1$, there are bounded Hankel operators without bounded symbols
(Theorem~2.1).  This indicates that $L^\infty $ symbols are not enough to
characterize bounded Hankel operators, and that \BMOr/ is the right class of
symbols in product domains [CS3].

For $d>1$ it is known [Am] that the positivity of the Pick matrix is
necessary but not sufficient for the existence of a solution of the
Pick problem.  Necessary and sufficient conditions involving Pick
matrices have been given by Agler for $d=2$ [Ag], and by Cole, Lewis
and Wermer for all $d>1$ [CLW].  However, their conditions are not
verifiable in practice, and the relation with Hankel operators is lost
in their approach.  In Section 3 we return to the consideration of
analogues to the Pick problem with \BMOr/-norm control initiated in
[CS3], and give necessary and sufficient conditions for the existence
of solutions of a coordinate-wise Pick problem in terms of either the
boundedness of a Hankel operator with symbol specified by the data, or
the positiveness of $d$ associated $n\times n$ Pick matrices.

In the case $d > 1$, all singular numbers of a Hankel operator are bounded
below by $d^{-1/2}$ times its norm (Theorem $B)$, so that all Hankel operators
of finite rank are zero [CS2].  This abrupt change from the
one-dimensional case is closely related to the failure of the Beurling
characterization of invariant subspaces to hold in the polydisk [AhC], and
shows that an AAK theory cannot be meaningful in $\T^d$, for $d > 1$.  To recover
the main features of the Nehari--AAK theory we need to introduce,
not only \BMOr/ symbols, but {\it sigma numbers\/} 
to replace the singular numbers, and a notion of
operators of {\it finite type}, to replace that of finite rank.

In Section 4 we rely on a version of Beurling's characterization in the
polydisk given in [CS4] to characterize the symbols of Hankel operators of
finite type in terms of tensor products of finite Blaschke products, and to
extend the AAK result mentioned above in terms of the sigma numbers of the
Hankel operators.

In Section 1 it is shown that, when passing from $\T$ to $\T^d$, for
$d>1$, the different equivalent characterizations of $\BMO/(\T)$
give rise to distinct spaces.  Similarly, the different
characterizations of Carleson measures in $\D$ give rise to different
notions in $\D^d$, for $d>1$.  One such characterization is that a
measure in $\D$ is Carleson if and only if a canonically associated
function is in \BMO/.  In Section 5 we extend this canonical
association to $d>1$, by defining Carleson--Nikolskii measures, and
proving that a measure is of this type if and only if a canonically
associated function is in \BMOr/. 

In the circle, the norms of Hankel operators of finite rank coincide
with the norms of multipliers acting in finite-dimensional model
subspaces, which in turn are determined by finite Pick matrices.  In
Section~6 we prove that the norms of Hankel operators of finite type
coincide with those of multipliers acting in corresponding model
subspaces, which now are not finite-dimensional but of bi-finite type,
like those appearing in Sections 4 and 5.  This significantly reduces
the number of steps required to verify norm boundedness.

\subsection*{Acknowledgements}

We want to thank Chandler Davis for extensive discussions with the
second author on duality, and Nikolai Nikolskii for helpful comments.
The last version of this paper was written while the second author was
a Research Professor of the Mathematical Sciences Research Institute
at Berkeley, and we are happy to acknowledge the hospitality received
there by both of us.

\subsection*{Basic Notations}

The following notations will be used throughout the paper.  For $d \geq  1$,
${\cal P} = {\cal P} (\T^d)$ is the class of trigonometric polynomials;
$\hat f$ represents the Fourier transform of $f$;
$$
\eqalign{
H^2(\T^d) &= \{f\in L^2(\T^d): \hat f(n_1,\dots, n_d) = 0 \hbox{ if $n_k < 0$
for some $k = 1,\dots,d$}\};\cr
H^2_{x_k} &= \{f\in L^2(\T^d) = \hat f(n) = 0 \hbox{ for $n_k < 0$}\};\cr
H^2(\T^d)^\perp  &= L^2(\T^d) \ominus  H^2(\T^d);}
$$
and the orthogonal projector $P: L^2 \rightarrow H^2(\T^d)$ is called
the analytic projector.

The $d$ shifts $S_k= S_{x_k}$, in $L^2(\T^d)$,
where $k = 1,\dots,d$, are defined by
$$
S_kf(x) = S_kf(x_1,\dots, x_d) := \exp (ix_k)f(x).
$$
In the case $d = 2$, we write $(x,y)$ for $(x_1,x_2)$ and $(m,n)$ for
$(n_1,n_2)$, and consider the subspaces of $H^2(\T^2)$ given by
$$
H^2_x = \{f\in L^2(\T^2): \hat f(m,n) = 0\hbox{ for }m < 0\}, \quad
H^2_{-x} = L^2 \ominus  H^2_x,
$$
and
$$
H^2_y = \{f\in L^2(\T^2): \hat f(m,n) = 0\hbox{ for }n < 0\}, \quad
H^2_{-y} = L^2 \ominus  H^2_y,
$$
as well as the projectors
$$
P_x: L^2 \rightarrow  H^2_x,  \quad P_{-x} := (I - P_x): L^2 \rightarrow
H^2_{-x}
$$
and
$$
P_y: L^2 \rightarrow  H^2_y, \quad P_{-y} := (I - P_y).
$$
The two shifts  $S_1$ and $S_2$ in $L^2(\T^2)$ satisfy
$$
S_1f(x,y) = e^{ix}f(x),\quad S_2f(x,y) = e^{iy}f(y).
$$
Observe that
$$
H^2(\T^2)^\perp  = H^2_{-x} + H^2_{-y} = H^2_{-x}
\dotplus
(H^2_{-y} \cap  H^2_x) = H^2_{-y}
\dotplus
(H^2_{-x} \cap  H^2_y).
$$

\section{Two Distinguished Subspaces of Product \BMO/}

An integrable function in $\T$ is of {\it bounded mean oscillation\/} if
\begin{equation}
\label{eq1.1}
{1\over |I|} \int  _I|f(x) - f_I|\,dx \leq  C\qquad \hbox{for all
intervals $I$},
\end{equation}
where $f_I = |I|^{-1} \int  _If(x)\,dx$.  The class \BMO/ of functions of
bounded mean oscillation is important in analysis.  It is closely related to
the Carleson measures and to the $A_p$ weights, as well as to bounded Hankel
operators.

A function $\phi $ is in $\BMO/ = \BMO/(\T)$ if and only if a
canonically associated measure $\mu $ in $\D$ is Carleson, and a
measure $\mu $ in $\D$ is Carleson if and only if $\Gamma 1\in \BMO/$
for a canonically associated Hankel operator $\Gamma $.  (See
definitions below.)  As there are different characterizations for the
elements of \BMO/ in $\T$, the same is true for Carleson measures in $\D$.

\BMO/ coincides with the space $L^\infty + HL^\infty $, where $H$ is
the Hilbert transform.  This characterization follows from Charles
Fefferman's famous duality result, asserting that \BMO/ is the (real)
dual of the Hardy space $H^1$.  Another way to prove
\begin{equation}
\BMO/ = L^\infty  + HL^\infty
\label{eq1.2}
\end{equation}
is through the characterizations of the weights $w$ for which $H$ is bounded in
$L^2(w)$ given by the $A_2$ condition and by the Helson--Szeg\H{o} theorem
[HS].

In passing from $\T$ to $\T^d$, for $d > 1$, the extension of the \BMO/ theory to
product domains presents various difficulties [ChF2].  S.-Y. Alice Chang
and Robert Fefferman were able to introduce a notion of {\it product $\BMO/ =
\BMO/(\T^d)$}, dual to the space $H^1(\D^d)$, and for which
an analogue of \eqref{1.2} is retained [ChF1].  In fact,
\begin{equation}
\eqalign{
\phi &\in \BMO/(\T^d)\iff{}\cr
\phi &= f_1 + H_{x_1}f_2 +\cdots+ H_{x_d}f_{d+1} + H_{x_1}H_{x_2}f_{d+2}
+\cdots+ H_{x_1}H_{x_2}\dots H_{x_d}f_{2^d}\cr
&\qquad\hbox{for $f_1,\dots, f_{2^d} \in L^\infty (\T^d)$,}}
\label{eq1.3}
\end{equation}
where $H_{x_j}$ is the Hilbert transform with respect to the variable
$x_j$, for $j = 1,\dots, d$, and \BMO/ is a complete normed space with
respect to
$$
\|\phi \|_{\BMO/} := \inf  \{\max_j
\|f_j\|_\infty :\hbox{ all decompositions } \eqref{1.3}\}.
$$
But for product \BMO/ the geometric characterizations
by mean oscillation and by
associated Carleson measures become considerably more complicated (they do not
correspond to bounded mean oscillation with respect to rectangles), and,
furthermore, the connections with weights and Hankel operators are lost.

In previous work ([CS1], [CS3]), we gave results in product spaces
analogous to those linking \BMO/ to weights and to Hankel operators in
one variable, in terms of classes of functions that are properly
contained in product \BMO/.  In this section we clarify the relation
of these classes with product \BMO/, give some of their basic
properties, and characterize their preduals.

\begin{definition}[small \BMO/]
\label{def1}
A function $\phi \in L^2(\T^d)$, for $d \geq  1$, is in $\bmo/(\T^d)$ if
there exist
$f_1,\dots, f_d,g_1,\dots, g_d \in L^\infty (\T^d)$ such
that
\begin{equation}
\phi  = f_1 + H_{x_1}g_1 =\cdots= f_d + H_{x_d}g_d
\label{eq1.4}
\end{equation}
and
$$
\|\phi \|_{\bmo/} := \inf  \{\max_{1\le j\le d}
\{\|f_j\|_\infty, \|g_j\|_\infty \}: \hbox{ all decompositions
\eqref{1.4}}\}.
$$
\end{definition}

Observe that $\|\phi \|_{\bmo/} = 0$ if and only if $\phi$ is
constant, and $\bmo//\C $ is
a complete normed space with respect  to
$\|\cdot \|_{\bmo/}$.

\begin{definition}[restricted \BMO/]
\label{def2}
A function $\phi \in L^2(\T^d)$, for $d \geq 1$, is in \BMOr/ if there
exist $\varphi _0,\varphi _1,\dots, \varphi _d \in L^\infty (\T^d)$
such that 
\begin{equation}
\label{eq1.5}
\left\{\eqalign{
&(I - P_{x_j})\phi = (I - P_{x_j})\varphi _j\quad\hbox{
for $j = 1,\dots,d$,\quad and}\cr
&P_{x_1}P_{x_2}\dots P_{x_d}\phi  = P_{x_1}P_{x_2}\dots P_{x_d}\varphi _0,}
\right.
\end{equation}
where $P_{x_j}: L^2 \rightarrow  H^2_{x_j}$ is the analytic projector in
$x_j$, for $j = 1,\dots, d$.
Moreover,
{\hfuzz=4pt
$$
\eqalign{
\|\phi \|_{\BMOr/} &:= \inf  \{\max_{0\le j\le d}
\|\varphi _j\|_\infty : \hbox{ all decompositions \eqref{1.5}}\}\cr
 &\phantom{:}=\max \{\max_{1\le j\le d}\{\inf
\{\|\phi  - h_{x_j}\|_\infty : h_{x_j} \in H^2_{x_j}\}\}, \,\inf
\{\|\phi  - h^\perp \|_\infty :  h^\perp  \in H^{2\perp }\}\}.}
$$}
\end{definition}

Observe that \BMOr/ is a complete normed space with respect to
$\|\cdot \|_{\BMOr/}$, and coincides with the space restricted \BMO/
introduced in [CS3].

The two definitions given above are
justified by the following results.

\begin{unnumberedtheorem}[Helson--Szeg\H{o} theorem in $\T^d$, for $d \geq 1$]
{\normalfont [CS1]}
A weight $0 \leq w\in L^1(\T^d)$ satisfies
$$
\int  _{\T^d}|Hf|^2w \leq  M^2\int  _{\T^d}|f|^2w\quad \hbox{
for all $f\in {\cal P}$},
$$
where $H=H_{x_1}\dots H_{x_d}$ is the product Hilbert transform,
if and only if $\phi  = \log w\in \bmo/(\T^d)$, with
$$
\phi  = u_1 + H_{x_1}v_1 =\cdots= u_d + H_{x_d}v_d
$$
for $u_1,\dots,u_d,v_1,\dots, v_d$ real-valued bounded functions in
$\T^d$ satisfying $\|u_j\|_\infty \leq C_M$ and $\|v_j\|_\infty \leq \pi /
2 - \epsilon _M$ for $j = 1,\dots, d$.
\end{unnumberedtheorem}

\begin{unnumberedtheorem}[Nehari theorem in $\T^d$, for $d \geq 1$]
{\normalfont [CS3]}
Let $\Gamma : H^2(\T^d)\cap\script P
\rightarrow  H^2(\T^2)^\perp $ be a Hankel operator. {\normalfont
(The definition of a Hankel operator is given in Section 2.)}
$\Gamma $ is bounded
if and only if there exists
$\phi \in \BMOr/$ satisfying $\Gamma f = (I -
P)(\phi f)$, for all $f\in H^2(\T^d)$, where $P: L^2 \rightarrow
H^2$ is the orthogonal projector, and $\|\phi \|_{\BMOr/} \leq  \|\Gamma \|
\leq  \sqrt{d} \|\phi \|_{\BMOr/}.$
\end{unnumberedtheorem}

Definitions 1 and 2 impose constraints on the functions in small and
restricted \BMO/, which follow immediately from the relation between
the analytic projector and the Hilbert transforms,
\begin{equation}
P_{x_j} = {1\over 2} (I + H_{x_j})\quad\hbox{for $j = 1,\dots,d$},
\label{eq1.7}
\end{equation}
and can be summarized as follows:

\begin{lemma}
\label{1.1}
\begin{enumerate}
\item[\normalfont (i)] For $\phi \in  \bmo/$ given by \eqref{1.4},
and for $j = 1,\dots,d$,    we have
\begin{equation}
\eqalign{
P_{x_j}\phi  &= P_{x_j}(f_j - ig_j)\cr
(I - P_{x_j}) \phi  &= (I - P_{x_j})(f_j + ig_j).}
\label{eq1.8}
\end{equation}
In particular,
\begin{equation}
P_{x_1}\dots P_{x_d}(f_1 - ig_1) =\cdots= P_{x_1}\dots P_{x_d}(f_d -
ig_d)
\label{eq1.9}
\end{equation}
and
\begin{equation}
\label{eq1.9a}
(I - P_{x_1})\dots (I - P_{x_d})(f_1 + ig_1) =\cdots= (I - P_{x_1})\dots (I -
P_{x_d})(f_d + ig_d).
\end{equation}

\item[\normalfont (ii)]
For $\phi \in \BMOr/$  {\it given by} \eqref{1.5},  we have
\begin{equation}
(I - P_{x_1})\dots (I - P_{x_d})\varphi _1 =\cdots= (I - P_{x_1})\dots (I -
P_{x_d})\varphi _d.
\label{eq1.9b}
\end{equation}
\end{enumerate}
\end{lemma}

Lemma \ref{1.1} implies that in order to define functions in \bmo/ or
\BMOr/ by $d$ or $d + 1$ bounded functions, respectively, those
bounded functions have to satisfy the constraints \eqref{1.9},
\eqref{1.9a} and \eqref{1.9b}.

The relation between small, restricted and product \BMO/s is the
following:

\begin{proposition}
\label{1.2}
The inclusions
$$
L^\infty (\T^d) \subset  \bmo/(\T^d) \subset  \BMOr/(\T^d) \subset  \BMO/(\T^d)
$$
are topological, and proper for $d > 1$.  For $d = 1$ we have
$\bmo/(\T) = \BMOr/(\T) = \BMO/(\T)$.
\end{proposition}

\begin{proof}
The topological inclusion $L^\infty  \subset  \bmo/$ is immediate
from Definition~\ref{def1}.
If $\phi \in $ \bmo/, by Lemma~\ref{1.1}(i), for $j = 1,\dots,d$,
we have $(I - P_{x_j})\phi  = (I - P_{x_j})(f_j + ig_j)$, and
$P_{x_1}\dots P_{x_d}\phi  = P_{x_1}\dots P_{x_d}(f_1 - ig_1) =\cdots= P_{x_1}\dots P_{x_d}(f_d - ig_d)$, with $f_j \pm  ig_j \in L^\infty (\T^d)$,
which means that condition \eqref{1.5} is satisfied and $\phi \in \BMOr/$  with
$\|\phi \|_{\bmo/} \geq  \|\phi \|_{\BMOr/}$.  It follows from
\eqref{1.5} and \eqref{1.7} that
$\phi \in \BMOr/$  implies $\phi \in \BMO/$, with $\|\phi \|_{\BMOr/}
\geq  \|\phi \|_{\BMO/}$.
To show that the inclusions are proper it is enough to
consider $d = 2$.

\begin{enumerate}
\item[(a)]  Example of $\phi \in \bmo/\backslash L^\infty $.  Fix
$v\in L^\infty (\T)$ such that $Hv \notin  L^\infty (\T)$, and let
$\phi (x,y) = Hv(x - y)$.  Then $\phi  = H_xg_1 = H_yg_2$, for $g_1,g_2
\in L^\infty (\T^2)$, defined by $g_1(x,y) = v(x - y)$,
$g_2(x,y) = -v(x - y)$, and $\phi \in \bmo/\backslash L^\infty (\T^2).$

\item[(b)]
Example of $\phi \in  \BMOr/\backslash \bmo/$.  Let
$\varphi _0(x,y) \equiv  0$, $\varphi _1(x,y) = v(x)h(y)$, $\varphi _2(x,y) =
h(x)v(y)$, where $v\in L^\infty (\T)$ is as in (a) and
$h\in H^\infty (\T)$ is not a constant.  Define $\phi \in  \BMOr/(\T^2)$
by condition \eqref{1.5}, that is,
$P_xP_y\phi  = P_xP_y\varphi _0$, $(I -
P_x)\phi  = (I - P_x)\varphi _1$, $(I - P_y)\phi  = (I - P_y)\varphi _2$, which
can be done since \eqref{1.9b} is satisfied:
$$
(I - P_x) (I - P_y)\varphi _1 = 0 = (I - P_x)(I - P_y)\varphi _2.
$$
If $\phi$ were in $\bmo/(\T^2)$, by
\eqref{1.8}, we would have $P_x\phi
= P_x(f_1 - ig_1)$ and $P_y\phi  = P_y(f_2 - ig_2)$, with $f_1 - ig_1$,
$f_2 - ig_2 \in L^\infty (\T^2)$.  But, in our case, $P_x\phi  =
P_xP_y\phi  + P_x(I - P_y)\phi  = P_xP_y\varphi _0
+ P_x(I - P_y)\varphi _2 = h(x)(I - P_y)v(y)$, with $(I - P_y)v(y) \notin
L^\infty (\T)$ by assumption.  Since
$P_x\phi  = P_x\varphi $ for some $\varphi \in L^\infty (\T^2)$, the
function $h(x)(I - P_y))v(y) = P_x\varphi (x,y)$ should be, for all $y$ fixed,
a function in $\BMO/(\T)$ satisfying $|(I - P_y)v(y)| \|h\|_{\BMO/(\T)} =
\|P_x\varphi (\cdot,y)\|_{\BMO/(\T)} \leq  c\|\varphi \|_\infty $.  Since
$\|h\|_{\BMO/} \neq  0$ and $(I - P_y)v(y) \notin  L^\infty (\T)$, this is a
contradiction.

\item[(c)]
Example of $\phi \in \BMO/\backslash $\BMOr/.  Given
$\psi _1,\psi _2 \in L^\infty (\T^2)$, take $\phi  = H_x\psi _1 + H_y\psi _2
\in \BMO/(\T^2)$.  In this case, $(I - P_x)\phi  = (I - P_x)(i\psi _1 +
i\psi _2) - (I - P_x)P_y(2i\psi _2)$, so, for $\phi $ to be in \BMOr/, by
\eqref{1.5}, $(I - P_x)P_y\psi _2$ should equal $(I - P_x)\varphi $, for
some $\varphi \in L^\infty (\T^2)$.  Taking
$\psi _2(x,y) = f(x)v(y)$, for $f, v\in L^\infty (\T)$, with $v$ as in
example (a), this means that, for all $y$ fixed, $|P_yv(y)| \|(I -
P_x)f\|_{\BMO/} = \|(I - P_x)P_y\varphi (\cdot,y)\|_{\BMO/(\T)} \leq
C\|\varphi \|_\infty $, which is a contradiction.
\end{enumerate}
\end{proof}

In what follows we limit the statements, as well as their proofs, to the case
$d = 2$, in order to simplify notations.  All results remain valid, with
obvious modifications, for $d > 1$.

We will write $\phi\in\BMO/_x(\T^2)$ if $\phi(\,\cdot,y)\in\BMO/(\T)$
for every $y$.  If, in addition, we have
$\sup_y\|\phi(\,\cdot,y)\|_{\BMO/}\le C$ for some constant $C$, we say
that $\phi\in\BMO/_x(\T^2)$ {\it with uniformly bounded norm}.  We
define $\BMO/_y(\T^2)$ similarly.

\begin{proposition}[Bounded mean oscillation on rectangles]
\label{1.2a}
The following conditions on a function $\phi$ are equivalent:
\emergencystretch5pt
\begin{enumerate}
\item[\normalfont(a)]  $\phi \in \bmo/(\T^2)$.

\item[\normalfont(b)]
For a constant $C > 0$ we have
\begin{equation}
{1\over |R|} \iint   _R|\phi (x,y) - \phi _R|\,dx\,dy \leq
C\quad \hbox{for all $R = I \times  J$},
\label{eq1.11}
\end{equation}
where $I,J \subset  \T$ are intervals and
$$
\phi _R = {1\over |R|} \iint   _R\phi (x,y)\,dx\,dy.
$$

\item[\normalfont(c)]
$\phi \in \BMO/_x(\T^2)$ with uniformly bounded norm and
$\phi \in \BMO/_y(\T^2)$ with uniformly bounded norm. 
\end{enumerate}
\end{proposition}

\begin{proof}
(b)$\implies$(c).  Condition \eqref{1.11}, of bounded mean oscillation
on rectangles, can be written as
$$
{1\over |I|} {I\over |J|} \int  _I \int  _J|\phi (x,y) - \phi _R|\,dx\,dy  =
{1\over |I|} \int  _IF(x,J)\,dx \leq  C
$$
for all intervals $I,J$.  This implies, for almost every $x\in $ I,
\begin{equation}
F(x,J) = {1\over |J|} \int  _J|\phi (x,y) - \phi _R|\,dy \leq  C\quad
\hbox{for all} J,
\end{equation}
which is to say that $\phi \in \BMO/_y$ with uniformly bounded norm.
Similarly, \eqref{1.11} implies $\phi \in \BMO/_x$ with uniformly
bounded norm.

\smallskip\noindent (c)$\iff$(a).  Obviously, $\phi \in \bmo/(\T^2)$
implies $\phi \in \BMO/_x(\T^2)$ and $\phi \in \BMO/_y(\T^2)$, both with uniformly
bounded norm.  Conversely, from (c) we have
$\phi=f_1+H_xg_1=f_2+H_yg_2$, 
for $f_1$ and $g_1$ bounded functions in $x$,
uniformly in $y$, and $f_2$ and $g_2$ bounded functions in $y$,
uniformly in $x$.  This means that
$f_1$, $f_2$, $g_1$ and $g_2$ are bounded
functions of both $x$ and in $y$, which is (a).

\smallskip\noindent
(c)$\implies$(b).
By \eqref{1.1}, the condition that $\phi \in \BMO/_x$ with uniformly
bounded norm is equivalent to
\begin{equation}
\label{eq1.12}
{1\over |I|} \int  _I|\phi (x,y) - \phi _I(y)|\,dx \leq  c\quad
\hbox{for all $I$, uniformly in $y$},
\end{equation}
and the condition that $\phi \in \BMO/_y$ with uniformly bounded norm
is equivalent to
\begin{equation}
{1\over |J|}\int  _J|\phi (x,y) - \phi _J(x)|\,dy \leq  c',\quad
\hbox{for all $J$, uniformly in $x$},
\label{eq1.12a}
\end{equation}
where
$$
\phi _I(y) = {1\over |I|} \int  _I \phi (x,y)\,dx,\qquad
\phi _J(x) = {1\over |J|} \int  _J \phi (x,y)\,dy,
$$
and $c,c'$ are positive constants.  From \eqref{1.12a} it follows, for
almost all $x\in I$, and $y\in J $, that
$$
|\phi (x,y) - \phi _J(x)| \leq  c' \quad\hbox{for all $J$},
$$
so that, for $R=I\times J$,
$$
\eqalign{|\phi _I(y) - \phi _R| &= \biggl|{1\over |I|}
\int  _I \phi (x,y)\,dx - {1\over |I|}
{1\over |J|} \int  _I \int  _J\phi (x,y)\,dx\,dy\biggr|   \cr
&= \bigl|{1\over |I|} \int  _I\phi (x,y)\,dx
- {1\over |I|} \int  _I \phi _J(x)\,dx\biggr|\cr
&\leq  {1\over |I|} \int  _I|\phi (x,y) - \phi _J(x)|\,dx \leq  c'.}
$$
Then, by \eqref{1.12},
$$
\displaylines{
{1\over |R|} \iint   _R |\phi (x,y) - \phi _R|\,dx\,dy \hfill\cr
\hfill\eqalign{
&= {1\over |J|}
\int  _J\biggl({1\over |I|} \int  _I |\phi (x,y) - \phi _R|\,dx\biggr)\,dy\cr
&\leq  {1\over |J|} \int  _J\biggl({1\over |I|} \int  _I|\phi (x,y) - \phi
_I(y)|\, dx + {1\over |I|} \int  _I |\phi _I(y) - \phi _R|\,dx\biggr)\,dy\cr
&\leq  c + c' = C,}}
$$
which is (b).
\end{proof}

The relation between \BMOr/, bounded Hankel operators in $\T^d$, for $d > 1$, and
Carleson measures will be treated in Section 4.  Now we consider duality
results.

In the one-dimensional case, $\BMO/(\T)$ is the dual of the (real) Hardy space
$$H^1(\T) := \{f\in L^1(\T): Hf\in L^1(\T)\},$$ or, equivalently,
is the space of functions $f$ such that $Pf$ and
$(I - P)f\in L^1(\T)$, where $f = P + (I - P)f$ is a canonical
decomposition of $f$ given by the analytic projector $P$.  In the
two-dimensional case, for each trigonometric polynomial $f\in {\cal P}
(\T^2)$, we consider three canonical decompositions of $f$, given in terms
of the analytic projectors $P_x: L^2 \rightarrow  H^2_x$ and $P_y: L^2
\rightarrow  H^2_y$, as well as of $P_{-x} := (I - P_x)$ and $P_{-y} := (I -
P_y)$:
\begin{eqnarray}
f &=& P_xf + P_{-x}f = P_yf + P_{-y}f;\label{eq(1)}\\
f &=& P_xf + P_{-x}P_yf + P_{-x}P_{-y}f = P_yf + P_{-y}P_xf +
P_{-y}P_{-x}f;\label{eq(2)}\\
f &=& P_xP_yf + P_xP_{-y}f + P_{-x}P_yf + P_{-x}P_{-y}f;\label{eq(3)}
\end{eqnarray}
and norm ${\cal P} (\T^2)$ with three different norms, all stronger
than the $L^1$ norm.  The completion of ${\cal P} (\T^2)$ with respect
to these three
norms gives rise to Banach spaces, denoted as follows:

\medskip\noindent(A)
The space
\begin{equation}
H^1_x(\T^2) + H^1_y(\T^2) := L^1(\T_y; H^1(\T_x)) + L^1(\T_x; H^1(\T_y)),
\label{eq(A)}
\end{equation}
whose elements are functions $f = f(x,y) = f_x(y)$, integrable in $y$, with
values in $H^1(\T_x)$, and
$f = f_y(x)$, integrable in $x$, with values in $H^1(\T_y)$; that
is, the closure of ${\cal P} (\T^2)$ in the norm
$$
[ f ] := \inf \{\|g\|_{H_x^1} + \|h\|_{H_y^1}: f = g + h\},
$$
where
$$
\|f\|_{H_x^1} := \|P_xf\|_1 + \|P_{-x}f\|_1\quad
\hbox{and}\quad \|f\|_{H_y^1} := \|P_yf\|_1 +
\|P_{-y}f\|_1
$$
correspond to partition \eqref{(1)}.

\medskip\noindent(B)
The space
\begin{equation}
{\cal H} ^1(\T^2) := {\cal H} ^1_x(\T^2) + {\cal H} ^2_y(\T^2),
\label{eq(B)}
\end{equation}
where ${\cal H} ^1_x$ and ${\cal H} ^1_y$ are the closures of ${\cal P} (\T^2)$
under the norms
corresponding to partition \eqref{(2)}, namely
$$
\eqalign{
\tri f\tri _{(x)} &:= \|P_xf\|_1 + (\|P_{-x}P_yf\|_1 + \|P_{-x}P_{-y}f\|_1),\cr
\tri f\tri _{(y)} &:= \|P_yf\|_1 + (\|P_{-y}P_xf\|_1 + \|P_{-y}P_{-x}f\|_1),}
$$
and ${\cal H} ^1$ is normed by
\begin{equation}
\tri f\tri  := \inf \{\tri g\tri _{(x)} + \tri h\tri _{(y)}: f = g + h\}.
\label{eq1.14a}
\end{equation}
Observe that, in particular, 
\begin{equation}
\tri f\tri \ge\|P_{-x}P_{-y}f\|_1 \quad\hbox{for all $f$}.
\label{eq1.14b}
\end{equation}

\medskip\noindent(C)
The space
\begin{equation}
H^1(\T^2) = H^1(\T_x; H^1(\T_y)) = H^1(\T_y; H^1(\T_x)),
\label{eq(C)}
\end{equation}
normed by
\begin{equation}
\|f\|_{H^1} := \|P_xP_yf\|_1 + \|P_xP_{-y}f\|_1 + \|P_{-x}P_yf\|_1 +
\|P_{-x}P_{-y}f\|_1,
\label{eq1.15}
\end{equation}
corresponding to partition \eqref{(3)}.

Observe that
$\cal P(\T^2)$ can be partitioned in more ways than those in
\eqref{(1)}--\eqref{(3)}.  For instance, a function $f$ can be written as
$f=P_xP_yf+(I-P_xP_y)f$, giving rise to the norm
$\|f\|:=\|P_xP_yf\|_1+\|(I-P_xP_y)f\|_1$.
Since the Hilbert transform, as well as the analytic projection, 
is unbounded in $L^1(\T)$, this norm $\|\cdot\|$ is not comparable to
those above, and in particular to $\tri\cdot\tri$.  

\begin{proposition}
\label{newprop}
For each $\epsilon>0$, there is an $f$ satisfying
$$
\epsilon\tri f\tri\ge\|f\|:=\|P_xP_yf\|_1+\|(I-P_xP_y)f\|_1,
$$
where $\tri\cdot\tri$ is defined in \eqref{1.14a}.
\end{proposition}

\begin{proof}
Consider $f(x,y)=\overline{u(x)}\,v(y)+v(x)\,\overline{u(y)}$, with $u$ an
inner function and $v\in L^1$ such that $\|v\|_1=1$, while
$\|w\|_1>\epsilon^{-1}$, for $w=(I-P)v$.  Since this $f$ satisfies
$P_xP_yf=0$ and
$P_{-x}P_{-y}f=\overline{u(x)}w(y)+w(x)\overline{u(y)}$, we have
$\|f\|=\|f\|_1\ge 2$, and, by \eqref{1.14b}, $\tri f\tri\ge
\|P_{-x}P_{-y}f\|_1$, which, after multiplying by $u(x)u(y)$, is equal
to
$$
\iint |u(x)w(x)+u(y)w(y)|\,dx\,dy>\int|w(x)|\,dx>\epsilon^{-1}.
$$
To justify the last inequality, it is enough to choose a test function
$G$ in the predual of $L^1(\T)$, such that $\int u(y)w(y)G(y)\,dy=0$,
since then 
$$
\hidewidth
\eqalignbot{
\iint |u(x)w(x)+u(y)w(y)|\,dx\,dy&\ge\sup_F
\iint (u(x)w(x)+u(y)w(y))F(x)G(y)\,dx\,dy\cr&=
\sup_F\int u(x)w(x)F(x)\,dx=\int |u(x)w(x)|\,dx.}\hidewidth
\proved
$$
\end{proof}

\begin{theorem}[Duality]
\label{1.3} The spaces defined in {\upshape \eqref{(A)}, \eqref{(B)}}
and\/ {\upshape\eqref{(C)}} are the preduals of the \BMO/s in $\T^2$.
More precisely:
\begin{enumerate}
\item[\normalfont(a)]
$\bmo/(\T^2)$ is the dual of  $H^1_x(\T^2) + H^1_y(\T^2)$.

\item[\normalfont(b)]
$\BMOr/(\T^2)$ is the dual of ${\cal H} ^1(\T^2) = {\cal H} ^1_x(\T^2) +
{\cal H}^1_y(\T^2)$.

\item[\normalfont(c)]
$\BMO/(\T^2)$ is the dual of  $H^1(\T^2)$ {\normalfont [ChF1]}.
\end{enumerate}
\end{theorem}

\begin{proof}
Note that, for any pair of functions $f$ and $\phi $ in the
variables $x$ and $y$ for which the integrals make sense,
$$
\int (P_{\pm x}f) (P_{\pm x}\phi )\,dx = \int (P_{\pm y}f)(P_{\pm
y}\phi )\,dy = 0,
$$
so that
\begin{equation}
\int   (P_{\pm x}f)(P_{\mp x}\phi )\,dx = \int   (P_{\pm x}f)\phi \,dx = \int
f(P_{\mp x}\phi )\,dx,
\label{eq1.16a}
\end{equation}
and similarly for $P_{\pm y}.$

\medskip\noindent
(a)  Let $\phi \in \bmo/(\T^2)$ 
and $f\in H^1_x + H^1_y$.  Then, by Lemma~\ref{1.1}(i) and \eqref{1.16a},
we have
$\int   f\phi \,dx\,dy = \int   g\phi \,dx\,dy + \int   h\phi \,dx\,dy$, where
$$
\eqalign{
\int   g\phi \,dx\,dy &= \int   (P_xg)\phi  + \int   (P_{-x}g)\phi  = \int
(P_xg)(P_{-x}\phi ) + \int   (P_{-x}g)(P_x\phi )\cr&
= \int   (P_xg)(P_{-x}\psi _1) + \int   (P_{-x}g)(P_x\varphi _1)\cr&
= \int   (P_xg)\psi _1 + \int   (P_{-x}g)\varphi _1,}
$$
for $\psi _1,\varphi _1 \in L^\infty (\T^2)$.
Similarly,
$$
\int   h\phi  = \int   (P_yh)\psi _2 + \int   (P_{-y}h)\varphi _2,
$$
for $\psi _2,\varphi _2 \in L^\infty (\T^2)$. Thus,
$$
\biggl|\int   g\phi \biggr| 
\leq  \|P_xg\|_1 \|\psi _1\|_\infty  + \|P_{-x}g\|_1
\|\varphi _1\|_\infty
$$
and
$$
\biggl|\int   h\phi \biggr| 
\leq  \|P_yh\|_1 \|\psi _2\|_\infty  + \|P_{-y}h\|_1
\|\varphi _2\|_\infty,
$$
which imply
$$
\biggl|\int   g\phi \biggr| 
  \leq  \|g\|_{H_x^1} \|\phi \|_{\bmo/}
\quad\hbox{and}\quad
\biggl|\int   h\phi \biggr|   \leq  \|h\|_{H_y^1} \|\phi \|_{\bmo/},
$$
and hence
$$
\biggl|\int   f\phi \biggr| \leq  [ f ]\|\phi \|_{\bmo/}.
$$

Conversely, if $\ell \in (H^1_x + H^1_y)^*$, then for every $f\in
{\cal P} (\T^2)$, we have $|\ell (f)| \leq C\|f\|_{H^1_x}$,
independently of $y$, and $|\ell
(f)| \leq C\|f\|_{H^1_y}$,
independently of $x$.  Since $(H^1_x)^* = \BMO/_x$ and $(H^1_y)^*
= \BMO/_y$, this implies, by Proposition~\ref{1.2a},
that $\ell $ can be given by a function in $\bmo/(\T^2)$.

\medskip\noindent
(b)  Let $\phi \in  \BMOr/(\T^2)$ be given by
$$
\phi  = P_xP_y\varphi _0 + P_xP_{-y}\varphi _2 + P_{-x}\varphi _1 =
P_xP_y\varphi _0 + P_{-x}P_y\varphi _1 + P_{-y}\varphi _2
$$
for $\varphi _0,\varphi _1,\varphi _2 \in L^\infty (\T^2)$, and let
$f\in {\cal H}^1_x + {\cal H}^1_y$.  Then, for $f = g + h$,
$g\in {\cal H}^1_x$, $h\in {\cal H}^1_y$, we have
$$
\int   g\phi  = \int   (P_{-x}P_{-y}g)\varphi _0 + \int
(P_{-x}P_yg)\varphi _2 + \int   (P_xg)\varphi _1
$$
and $|\int g\phi | \leq \tri g\tri _{(x)} \|\phi \|_{\BMOr/}$.
Similarly, $|\int h\phi | \leq \tri h\tri _{(y)} \|\phi \|_{\BMOr/}$;
hence 
$$
\left|\int f\phi \right| \leq \tri f\tri\, \|\phi \|_{\BMOr/}.
$$

Conversely, if $\ell $ is continuous on  ${\cal H} ^1(\T^2)$, it is continuous
on ${\cal H} ^1_x$ and on ${\cal H}^1_y$.  In particular, for all
$f\in {\cal P} (\T^2)$,
$$
|\ell  (f)| \leq  C\tri f\tri _{(x)} = C\|P_xf\|_1 + C(\|P_yP_{-x}f\|_1 +
\|P_{-y}P_{-x}f\|_1).
$$
In particular, if $f = P_xf$, we have $|\ell  (f)| \leq  C\|f\|_1$, and
there exists $F\in L^\infty $ such that
$$
\ell  (P_xf) = \int   F(P_xf)\,dx\,dy = \int   (P_{-x}F)f \,dx\,dy,
$$
by \eqref{1.16a}.
Similarly, if $f = P_{-x}f$,
$$
|\ell  (f)| \leq  C(\|P_yf\|_1 + \|P_{-y}f\|_1) = C\|f\|_{H^1_y},
$$
and there exists $G\in \BMO/_y$ such that
$$
\ell  (P_{-x}f) = \int   G(P_{-x}f)\,dx\,dy = \int   (P_xG)f \,dx\,dy.
$$
Then
\begin{equation}
\ell  (f) = \ell  (P_xf) + \ell  (P_{-x}f) = \int   (P_{-x}F + P_xG)f
\,dx\,dy
\label{eq1.17}
\end{equation}
for $F\in L^\infty $ and $G\in \BMO/_y$.  Similarly,
\begin{equation}
\ell  (f) = \ell  (P_yf) + \ell  (P_{-y}f) = \int   (P_{-y} F' + P_yG')f
\,dx\,dy
\label{eq1.17a}
\end{equation}
for $F'\in L^\infty $ and $G'\in \BMO/_x$.

Since, by \eqref{1.17} and \eqref{1.17a}, the two functions
representing $l$ coincide as functionals on all $f\in {\cal P} (\T^2)$,
we conclude that
$$
P_{-x}F + P_xG = P_{-y}F' + P_yG' = \phi .
$$

Now observe that $G\in \BMO/_y$ and $G'\in \BMO/_x$ imply that $P_xG =
P_x\varphi $ and $P_yG' = P_y\varphi '$ for some $\varphi,\varphi
'\in L^\infty (\T^2)$.  Therefore, the function $\phi $ satisfies
$$
P_xP_y\phi  = P_xP_y\varphi  = P_xP_y\varphi '\quad\hbox{ and }\quad
P_{-x}\phi  = P_{-x}F,\quad P_{-y}\phi  = P_{-y}F'
$$
for $\varphi,\varphi ', F, F'\in L^\infty (\T^2)$, which means, by
definition, 
that $\phi \in  \BMOr/(\T^2)$.

\medskip \noindent
(c)  Let $\phi \in \BMO/(\T^2) = L^\infty  + H_xL^\infty  +
H_yL^\infty  + H_xH_yL^\infty$ (see  [ChF1]) and $f\in H^1$.  Writing
$$
\phi  = P_xP_y\psi _1 + P_xP_{-y}\psi _2 + P_{-x}P_y\psi _3 +
P_{-x}P_{-y}\psi _4,
$$
for $\psi _1,\psi _2,\psi _3,\psi _4 \in L^\infty $, we get
$$
\displaylines{
\biggl|\int   f\phi \biggr| = \biggl|\int
(P_xP_yf)(P_{-x}P_{-y}\psi _4) + \int
(P_xP_{-y}f)(P_{-x}P_y\psi _3) \hfill\cr
\hfill {}+ \int   (P_{-x}P_yf)(P_xP_{-y}\psi _2)
+ \int   (P_{-x}P_{-y}f)(P_xP_y\psi _1)\biggr| \cr
\phantom{\biggl|\int   f\phi \biggr|}
\leq  \|P_xP_yf\|_1 \|\psi _4\|_\infty  + \|P_xP_{-y}f\|_1 \|\psi _3\|_\infty
\hfill\cr\hfill{}+ \|P_{-x}P_yf\|_1 \|\psi _2\|_\infty  
+ \|P_{-x}P_{-y}f\|_1 \|\psi _1\|_\infty ,}
$$
which implies $\bigl|\int f\phi\bigr| \leq  \|f\|_{H^1} \|\phi \|_{\BMO/}$.

Conversely, if $\ell \in  (H^1)^*$, the usual duality
argument shows that $\ell$ is given by a function in
$L^\infty+H_x L^\infty+H_y L^\infty+H_xH_y L^\infty=
\BMO/(\T^2)$, as in [ChF1].
\end{proof}

A more detailed study of \bmo/ and \BMOr/ in $\T^d$, for $d > 1$, including
their atomic decompositions and their associated Carleson measures,
will be the object of a future paper.

\section{Big Hankel Operators and Their \BMOr/ Symbols}

We consider operators $\Gamma : {\cal P}  \cap  H^2(\T^d) \rightarrow
H^2(\T^d)^\bot $, for $d \geq  1$. Such operators
$\Gamma $ are called bounded if
$\sup \|\Gamma f\|_2/\|f\|_2 =: \|\Gamma \| < \infty $.  A bounded $\Gamma $
has a unique bounded extension, $\Gamma : H^2 \rightarrow  H^{2\bot }$.  It is
easy to check that, for every $\Gamma : {\cal P}  \cap  H^2 \rightarrow
H^{2\bot }$, the following conditions are equivalent$:$

\begin{enumerate} 
\item[(a)] 
$\<\Gamma S_kf$, $g\> = \<\Gamma f, S_{-k}g\>$ for $k =
1,\dots,d$ and all $f\in {\cal P} \cap H^2$ and $g\in H^{2\bot
}$;

\item[(b)] 
$\Gamma S_kf = (I - P)S_k\Gamma f$ for $k = 1,\dots, d$, where
$P: L^2 \rightarrow H^2$ is the orthoprojector;

\item[(c)] 
There exists $\phi \in L^2(\T^2) \ni \Gamma = \Gamma _\phi $, that
is, $\Gamma f = (I - P)\phi f$ for all $f\in H^2(\T^2)$;

\item[(d)] 
$\Gamma  = \Gamma _{\phi _-}$ for $\phi _- = \Gamma 1\in H^{2\bot }$.
\end{enumerate}

If (a)--(d) are verified, $\Gamma $ is called a {\it big Hankel
operator}, and $\phi $ as in (c) is called a {\it symbol\/} of $\Gamma $.
Since
$$
\Gamma _\phi  = \Gamma _\psi  \iff   \phi  - \psi  = h\in H^2,
$$
we see that if $\phi $ is a symbol for $\Gamma $ so are all $\phi +
h$, for $h\in H^2$.  Moreover, among all symbols, there is a unique
one in $H^{2\bot }$, which is $\Gamma 1$.  In what follows, (big)
Hankel operators will be referred to as {\it Hankel}.

If $\varphi \in L^\infty $, then $\Gamma _\varphi $ is a bounded
operator, with $\|\Gamma _\varphi  \leq  \|\varphi \|_\infty $.  In the
one-dimensional case, the Nehari Theorem gives the converse:  A Hankel
operator $\Gamma $ is bounded if and only if $\exists \varphi \in L^\infty $ with
$\Gamma _\varphi  = \Gamma $, if and only if $\exists \varphi \in L^\infty $ with
$\Gamma _\varphi  = \Gamma $ and $\|\varphi \|_\infty  = \|\Gamma \|$,  and if and only if
$\Gamma 1\in  \BMO/$.  Also,
$\|\Gamma _\varphi \| = \dist_{L^\infty }(\varphi,H^\infty )$.  Since
$\phi \in \BMO/(\T)$ implies $\phi  = \varphi  + h$, for
$\varphi \in L^\infty $ and, $h\in H^2$, we have
$$
\phi \in \BMO/ \implies \Gamma _\phi = \Gamma _\varphi \hbox{ is
bounded, and } \|\Gamma _\phi \| \leq \|\phi \|_{\BMO/}.
$$

Thus, in the one-dimensional case, $\BMO/(\T)$ appears as an essential
feature both in the weighted norm inequalities for the Hilbert
transform, and in the boundedness of the Hankel operators.  In [ACS]
it was shown that the basic properties of $\BMO/(\T)$ can be deduced in
a unified way from a Generalized Bochner Theorem (GBT) that is
equivalent to the Nehari theorem in $H^2(\T;\mu )$, and which unifies
the results of Nehari and Helson--Szeg\H{o}.  An abstract version of
this GBT led to a version of the Nehari theorem in $\T^d$ in terms of
\BMOr/, and to an extension of the Helson--Szeg\H{o} theorem in terms
of $\bmo/(\T^d)$.  Since $\bmo/ \neq L^\infty $ and \BMOr/ $\neq \BMO/$
for $d > 1$, this underlines the importance of these two subspaces of
product \BMO/.  

\medskip 

\noindent Here we will base our considerations on the
two-dimensional version of Nehari theorem:

{\def\thelemma{A}\addtocounter{lemma}{-1}
\begin{theorem}{\normalfont [CS1], [CS2]}
For every $($big$)$ Hankel operator $\Gamma : {\cal P} \cap H^2(\T^2)
\rightarrow H^2(\T^2)^\bot $, the following conditions are equivalent:
\begin{enumerate}
\item[\normalfont(a)]  
$\Gamma $  is bounded.
\item[\normalfont(b)]  
There exist  $\varphi _1,\varphi _2 \in L^\infty (\T^2)$, with
$\max \{\|\varphi _1\|_\infty, \|\varphi _2\|_\infty \} \leq  \|\Gamma \|$,
and such that
\begin{equation}
P_{-x}\Gamma  = P_{-x}\Gamma _{\varphi _1},\quad P_{-y}\Gamma  =
P_{-y}\Gamma _{\varphi _2}.
\label{eq2.4}
\end{equation}
\item[\normalfont(c)]  
There exists
$\phi \in  \BMOr/ \cap  H^{2\perp }$ with
$\Gamma  = \Gamma _\phi $ and $\|\phi \|_{\BMOr/} \leq  \|\Gamma \| \leq
\sqrt{2}\|\phi \|_{\BMOr/}.$
\item[\normalfont(d)]  
$\Gamma 1\in  \BMOr/$.
\end{enumerate}
\end{theorem}}

Remark that \eqref{2.4} implies, for $\varphi _1$ and $\varphi _2$ as in (b),
that
\begin{equation}
\Gamma  = P_{-x}\Gamma _{\varphi _1} + P_xP_{-y}\Gamma _{\varphi _2} =
P_{-y}\Gamma _{\varphi _2} + P_{-x}P_y\Gamma _{\varphi _1}.
\label{eq2.4a}
\end{equation}

\begin{unnumberedcorollary} 
For every $\phi \in  \BMOr/$, $\Gamma _\phi $ is bounded and
\begin{equation}
\label{eq2.5}
\|\Gamma _\phi \| = \max  \{\dist_{\BMOr/}(\phi, \BMOr/ \cap  H^2_x),
\dist_{\BMOr/}(\phi, \BMOr/ \cap  H^2_y)\}.
\end{equation}
\end{unnumberedcorollary} 

From Theorem A and the fact that $\phi \in \BMOr/$ implies $\phi _- =
(I - P)\phi \in \BMOr/ \cap H^{2\bot }$, with $\|\phi _-\|_{\BMOr/}
\leq \|\phi \|_{\BMOr/}$, it follows that, for all $\phi \in \BMOr/$,
the operator $\Gamma _\phi$ is bounded and satisfies $\|\Gamma _\phi
\| \leq \sqrt{2} \|\phi \|_{\BMOr/}$.

Then, $\phi  \mapsto  \Gamma _\phi $ is a surjective map from \BMOr/ onto the
space ${\cal G}$  of the bounded Hankel operators, whose restriction to $\BMOr/
\cap  H^{2\bot }$ is a bijection.  The symbols $\phi \in  \BMOr/$, a
proper subspace of product \BMO/, are thus enough for the theory of big Hankel
operators.  The duality theorem for \BMOr/ leads to the following theorem, 
which
highlights that the symbols in $L^\infty $ are not enough, so that equivalence
(c) in Theorem A can be considered sharp.

The map $\varphi  \mapsto  \Gamma _\varphi $ from $L^\infty (\T^2)$ to the space
${\cal G}$  has kernel $H^\infty (\T^2)$ and induces an injective map from
$L^\infty /H^\infty $ into ${\cal G}$ .  If this map were also surjective, by
the Banach open mapping theorem, there would be a constant $K > 0$ such that,
for each $\Gamma  \in {\cal G}$, there would be a
$\varphi \in L^\infty $ with
\begin{equation}
\Gamma _\varphi  = \Gamma \quad\hbox{  and }\quad  \|\varphi \|_\infty  \leq
K\|\Gamma \|.
\label{eq2.6}
\end{equation}

If \eqref{2.6} held, the space $L^\infty+PL^\infty$, the ``natural''
extension of $\BMO/(\T)$ to $d>1$, would indeed coincide with \BMOr/.
We see that this is not the case by showing that the map $\varphi
\mapsto \Gamma $ is not surjective from $L^\infty /H^\infty $ to
${\cal G}$ .

\begin{theorem}  
\label{2.1}
There are bounded big Hankel operators from $H^2(\T^2)$
to $H^2(\T^2)^\bot $ that have no bounded symbol.
\end{theorem}  

\begin{proof}  
If the map $\varphi \mapsto \Gamma $ from $L^\infty /H^\infty $ to
${\cal G}$ were surjective, there would exist a $K > 0$ for which
\eqref{2.6} would be satisfied.  If a pair $\varphi _1,\varphi _2 \in
L^\infty (\T^2)$, with $\|\varphi _1\|_\infty \leq 1$ and $\|\varphi
_2\|_\infty \leq 1$, coincide as functionals on $H^2(\T^2)$, by
\eqref{2.4} it defines a bounded Hankel operator $\Gamma $, with
$\|\Gamma \| \leq \sqrt{2}$, so that there would be a $\varphi \in
L^\infty $ with $\Gamma = \Gamma _\varphi $ and $\|\varphi \|_\infty
\leq \sqrt{2}K$.  Now given any trigonometric polynomials $p_0 \in
H^2_{-x,-y} := H^2_{-x} \cap H^2_{-y}$, $p_1 \in H^2_{-x,y} :=
H^2_{-x} \cap H^2_y$, and $p_2 \in H^2_{x,-y} := H^2_x \cap H^2_{-y}$,
we have by \eqref{2.4a}:
$$
\eqalign{
\<\Gamma 1, \,p_0 + p_1 + p_2\> &= \<P_{-x}\Gamma _{\varphi _1}1, \,p_0
+ p_1\> + \<P_{-y}P_x\Gamma _{\varphi _2}1, p_2\>\cr
&=\<\Gamma _{\varphi _1}1, \,p_0 + p_1\> + \<\Gamma _{\varphi _2}1, p_2\>
= \int   \varphi _1(\bar p_0 + \bar p_1) + \int   \varphi _2\bar p_2.}
$$
On the other hand, if \eqref{2.6} holds, we have
$$
\eqalign{
|\<\Gamma 1, p_0 + p_1 + p_2\>| &=
|\<\Gamma _\varphi 1, p_0 + p_1 + p_2\>| = |\int   \varphi (\bar p_0 +
\bar p_1 + \bar p_2)|\cr
&\leq  \|\varphi \|_\infty  \| p_0 +  p_1 +  p_2\|_1 \leq
\sqrt{2}K\| p_0 +  p_1 +  p_2\|_1,}
$$
so that
\begin{equation}
\label{eq2.7}
\biggl|\int   \varphi _1( \bar p_0 +  \bar p_1) + \int   \varphi _2 \bar p_2\biggr|
\leq  \sqrt{2} K\| p_0 +  p_1 +  p_2\|_1.
\end{equation}

We will now show that \eqref{2.7} 
leads to a contradiction.  In fact, to give any pair
$\varphi _1,\varphi _2$ as above is the same as to give a $\phi _- \in 
\BMOr/ \cap  H^{2\bot }$, with $P_{-x}\phi _- = P_{-x}\varphi _1$,
$P_{-y}\phi _- = P_{-y}\varphi _2$, and
$P_xP_y\phi _- = 0$.  Then \eqref{2.7} can be rewritten as
\begin{equation}
\label{eq2.7a}
\biggl|\int   
\phi _-( \bar p_0 +  \bar p_1 +  \bar p_2)\biggr| \leq  \sqrt{2}K
\| p_0 +  p_1 +  p_2\|_1
\end{equation}
for all trigonometric polynomials $p_0 + p_1 + p_2 \in H^2_{-x,-y} +
H^2_{-x,y} + H^2_{x,-y} = H^{2\bot }$.

Now, any $\phi \in  \BMOr/$ can be written as $\phi  = P_xP_y\phi  + (I -
P_xP_y)\phi  = P_xP_y\varphi _0 + \phi _-$ for some $\varphi _0
\in L^\infty $ and $\phi _- \in  \BMOr/ \cap  H^{2\bot }$.  Thus,
for every trigonometric polynomial
$f = \bar p + \bar p_0 + \bar p_1 + \bar p_2$, 
with $p\in H^2(\T^2)$, and $\phi \in  \BMOr/$, \eqref{2.7a} yields
\begin{equation}
\thickmuskip=.5 \thickmuskip
\medmuskip=.5 \medmuskip
\!\!
\eqaligntop{\Bigl|\int \phi f \Bigr| & = \Bigl|\int \varphi _0 \bar p
+ \int \phi _-( \bar p_0 + \bar p_1 + \bar p_2)\Bigr| \leq \|\varphi
_0\|_\infty \| p\|_1 + \sqrt{2} K\| p_0 +  p_1 + 
p_2\|_1\cr& \leq (\|\phi \|_{\BMOr/} + \sqrt{2} K) (\|P_xP_y f\|_1 +
\|(I-P_xP_y) f\|_1),} \!\!
\label{eq2.8}
\end{equation}
where $p = P_xP_yf$ and $p_0 + p_1 + p_2 = (I - P_xP_y)f$.

But since the Hilbert transforms, as well as the analytic projections,
are unbounded in $L^1$, the norms $\|P_xP_y f\|_1 + \|(I-P_xP_y)
f\|_1$ and $\tri f\tri$ are not comparable 
(see Proposition~\ref{newprop}), and there exists for every $\epsilon>0$ an $f\in
{\cal P} (\T^2)$ such that
$$
\|P_xP_y f \|_1 + \|(I - P_xP_y) f \|_1 < {\epsilon \over 1 + \sqrt{2}K}
\tri f\tri .
$$  
By Hahn--Banach and the duality of \BMOr/, there exists $\phi \in $
\BMOr/ such that $\|\phi \|_{\BMOr/} \leq 1$ and $\int \phi  f =
\tri f\tri $, so that \eqref{2.8} implies $\int \phi f < \epsilon \int
\phi f$, which is a contradiction.
\end{proof}

An important open question is whether for every
$\varphi \in L^\infty (\T^2)$ there is another
$\psi \in L^\infty (\T^2)$ such that $\Gamma _\psi  = \Gamma _\varphi $
and $\|\psi \|_\infty  \leq  K\|\Gamma _\varphi \|$, where $K$ is a universal
constant and $\Gamma _\varphi $ is the big Hankel operator defined by
$\T_\varphi f = (I - P)\varphi f$.  Some geometric properties of
$H^\infty (\T^2)$ and \BMOr/ make highly improbable a positive answer to this
question, which will be considered elsewhere.

\section{Interpolation Problems in the Polydisk, Hankel Operators and
Pick Matrices}

A basic interpolation problem in $\D^d$, for $d\ge1$, is the {\it Pick
problem\/}: Given $z_1,\dots,z_n\in \D^d$ and
$\lambda_1,\dots,\lambda_n\in \C$, find a necessary and sufficient
condition for the existence of an analytic function $F$ on $\D^d$
satisfying $F(z_k)=\lambda_k$, for $k=1,\dots,n$, with
$\|F\|_\infty\le 1$.  This problem can be reformulated in a way that
is slightly more general only for $d>1$, as follows: Given
$z_1,\dots,z_n\in \D^d$ and $G\in H^\infty(\D^d)$, find an analytic $F$
satisfying $F(z_k)=G(z_k)$ for $k=1,\dots,n$, and $\|F\|_\infty\le 1$.

In the case of $d=1$, the problem was solved by G. Pick in 1916, in
terms of the positivity of an associated $n\times n$ matrix given by
the data.  Another solution has been given in terms of the boundedness
of an associated Hankel operator given by the data.

\begin{theorem}[Pick]
Given $z_1,\dots,z_n\in \D$ and
$G\in H^\infty(\D)$, the following assertions are equivalent:
\begin{enumerate}
\item[(i)]
The Pick matrix
\begin{equation}
((1-G(z_j)\overline{G(z_k)})(1-z_j\overline{z_k})^{-1})_{j,k=1,\dots,n}
\label{eqnew3.1}
\end{equation}
is positive definite.
\item[(ii)]
The Hankel operator $\Gamma_\phi$ with symbol $\phi=\bar b G$, where $b$
is the Blaschke product with simple zeros at $z_1,\dots,z_n$, is
bounded, and $\|\Gamma_\phi\|\le1$.
\item[(iii)]
The Pick problem has a solution.
\end{enumerate}
\end{theorem}

The equivalence of (i) and (iii) was proved in [P], and that of (ii)
and (iii) can be obtained as a corollary of the Nehari theorem (see,
for instance, [Ni]).

For $d>1$ it is known [Am] that the positivity of the Pick matrix
analogous to \eqref{new3.1} is necessary
but not sufficient for the existence of a solution to the Pick
problem.  Necessary and sufficient conditions 
involving Pick matrices have been given by Agler for $d=2$ [Ag], and
by Cole, Lewis and Wermer for all $d>1$ [CLW].  However, their
conditions are not verifiable in practice.  Moreover, in their
approach the relation with Hankel operators is lost.

As the Nehari theorem for $d>1$ can be recovered by replacing the
$L^\infty$ norm by the $\BMOr/$ norm, considering the Pick problem
with \BMOr/-norm control allowed us in [CS3] to retain the
relation with Hankel operators (within a constant $\sqrt d$), but not
a Pick condition.

The following result, which also reduces to the Pick theorem when
$d=1$, gives necessary and sufficient
conditions for the existence of solutions of a coordinate-wise Pick
problem in terms of either the boundedness of a Hankel operator with
symbol specified by the data, or the positivity of $d$ associated
$n\times n$ Pick matrices.  We state it here only for $d=2$, but it
holds for all $d>1$, with obvious changes.

\begin{theorem}
\label{new3.1}
Given $(z_1,w_1),\dots,(z_n,w_n)\in \D^2$ and $G\in H^\infty(\D^2)$, let
$\bee=b_1\otimes b_2$, where $b_1$ and $b_2$ are finite
one-dimensional Blaschke products with simple zeros at $z_1,\dots,z_n$
and $w_1,\dots,w_n$, respectively.  The following assertions are
equivalent (up to a constant $\sqrt 2$):
\begin{enumerate}
\item[(i)]
The Pick matrices
\begin{equation}
((1-G(z_j,y)\overline{G(z_k,y)})(1-z_j\overline{z_k})^{-1})_{j,k=1,\dots,n}
\label{eqnew3.2}
\end{equation}
and
\begin{equation}
((1-G(x,w_j)\overline{G(x,w_k)})(1-w_j\overline{w_k})^{-1})_{j,k=1,\dots,n}
\label{eqnew3.2a}
\end{equation}
are positive definite for every $y\in\T$ and every $x\in\T$, respectively.
\item[(ii)]
The Hankel operator $\Gamma _\phi$, for $\phi=\bar\bee G$, is bounded,
with $\|\Gamma_\phi\|\le 1$.
\item[(iii)]
There is a function $F\in H^2(\D^2)$ satisfying
$F(z_k,w)=G(z_k,w)$ and $F(z,w_k)=G(z,w_k)$, for $k=1,\dots,n$, with
$\|\bar \bee F\|_{\BMOr/}\le 1$.
\item[(iv)]
There exist two bounded functions on $\D^2$, $F_1$ analytic in $z$ and
$F_2$ analytic in $w$, satisfying $F_1(z_k,y)=G(z_k,y)$ and
$F_2(x,w_k)=G(x,w_k)$, for $k=1,\dots,n$, with $\|F_1\|_\infty\le 1$
and $\|F_2\|_\infty\le 1$.
\end{enumerate}

More precisely, {\upshape (ii)} implies {\upshape (iii)} and {\upshape
(iv)}, and either {\upshape (iii)} or {\upshape (iv)} imply {\upshape
(ii)} with $\|\Gamma_\phi\|\le\sqrt2$, while {\upshape (i)} is
equivalent to {\upshape (iv).  (Compare [BH].)}
\end{theorem}

\begin{remark}
Observe that the loss of the $L^\infty$ norm in condition (iii) is
compensated by the strengthening of the interpolation condition to
each variable independently.
\end{remark}

\begin{proof}
(ii)$\implies$(iii).
If $\|\Gamma_\phi\|\le1$, by Theorem A there exists $\psi\in\BMOr/$
with $\|\psi\|_\BMOr/\le 1$, such that $\Gamma_\phi=\Gamma_\psi$.
Therefore $\psi-\phi=\psi-\bar bG\in H^2(\D^2)$.  Setting $F=b\psi=
G+b_1b_2h\in H^2(\D^2)$ it is immediate that $F$ satisfies all the
conditions of (iii).

\medskip\noindent
(iii)$\implies$(ii).
If $F$ satisfies the interpolation conditions (iii), then
$$F-G=b_1(z)h_1(z,w)=b_2(w)h_2(z,w)$$ for $h_1,h_2\in H^2(\D^2)$.  This
implies that there is an $h\in H^2(\D^2)$ such that
$F-G=b_1(z)b_2(w)h(z,w)$.  Thus, setting $\psi=\bar b F$, we have
$\|\psi\|_{\BMOr/}\le1$, and, by Theorem A, $\|\Gamma_\psi\|\le \sqrt2$.
But $\Gamma_\psi=\Gamma_\phi$, since $\psi-\phi=\bar b(F-G)=h$, so 
$\|\Gamma_\phi\|\le \sqrt2$.

\medskip\noindent (ii)$\implies$(iv).  By Theorem A(b), the condition
$\|\Gamma_\phi\|\le1$ implies that there exist $\varphi_1,\varphi_2\in
L^\infty(\T^2)$ with $\|\varphi_j\|_\infty\le 1$ for $j=1,2$, such that
$\phi=\bar b G=\phi_1+h_x=\phi_2 +h_y$ for $h_x\in H^2_x$ and $h_y\in
H^2_y$.  The functions $F_1=b_1\varphi_1$ and $F_2=b_2\varphi_2$ satisfy
$\|F_j\|_\infty\le1$ for $j=1,2$, as well as $F_1(z_k,w)=G(z_k,w)$ and
$F_2(z,w_k)=G(z,w_k)$.  Moreover,
$F_1=\overline{b_2(w)}G(z,w)+b_1(z)h_x(z,w)$ is analytic in $z$, and
$F_2=\overline{b_1(z)}G(z,w)+b_2(w)h_y(z,w)$ is analytic in $w$.

\medskip\noindent (iv)$\implies$(ii).  
By the interpolation conditions satisfied by $F_1$ and $F_2$, for each
$y\in\T$ we have $G(z,y)-F_1(z,y)=b_1(z)h_x(z,y)$ for $h_x\in H_x^2$,
and, for each $x\in\T$, we have
$G(x,w)-F_2(x,w)=b_2(w)h_y(x,w)$ for $h_y\in H_y^2$.  Hence, setting
$\varphi_j=\bar bF_j$ for $j=1,2$, we have $\|\varphi_j\|_\infty\le1$, and
$$
\eqalign{
\phi-\varphi_1&=\bar b(G-F_1)=\overline{b_2}h_x\in H_x^2,\cr
\phi-\varphi_2&=\bar b(G-F_2)=\overline{b_1}h_y\in H_y^2.\cr}
$$
Again by Theorem A, this implies $\|\Gamma_\phi\|\le \sqrt 2$.

\medskip\noindent (i)$\iff$(iv).  Apply Pick's Theorem to the
one-variable functions $G(z,\cdot)$ and $G(\cdot, w)$ separately.
Then the two solutions $F_1(z,\cdot)$ and $F_2(\cdot, w)$
satisfy $\|F_1(\cdot,w)\|_\infty\le 1$ for all $w\in\D$, and
$\|F_2(z,\cdot)\|_\infty\le 1$ for all $z\in\D$, so that
$\|F_1\|_\infty\le 1$ and $\|F_2\|_\infty\le 1$.  Conversely, observe
that (iv) implies the analogues of \eqref{new3.2} and
\eqref{new3.2a} with $w\in\D$ instead of $y\in\T$, and $z\in\D$
instead of $x\in\T$, respectively, which is equivalent to
\eqref{new3.2} and
\eqref{new3.2a} by the analyticity  of $G$ in both variables.
\end{proof}

Given $(z_1,w_1),\dots,(z_n,w_n)\in \D^2$ and
$\lambda_1,\dots,\lambda_n$ in $\C$, let $b_1$ and $b_2$ be the
corresponding one-dimensional Blaschke products.  Writing
$\lambda_k=\lambda_k' \lambda_k''$, for $k=1,\dots,n$, set
$$
G_1(z)=\sum_{k=1}^n{b_1(z)\over b'_1(z_k)}{\lambda_k'\over z-z_k}
\qquad\hbox{and}\qquad
G_2(w)=\sum_{k=1}^n{b_2(w)\over b'_2(w_k)}{\lambda_k''\over w-w_k},
$$
so that $G_1(z_k)=\lambda_k'$, $G_2(w_k)=\lambda_k''$, and
$G(z_k,w_k)=\lambda_k$ for $G(z,w)=G_1(z) G_2(w)$.  For such $G_1$ and
$G_2$, or any others satisfying the interpolating conditions, we have:

\begin{corollary}
Given $(z_1,w_1),\dots,(z_n,w_n)\in \D^2$ and
$\lambda_1,\dots,\lambda_n$ in $\C$, there exists a function $F\in
H^2(\D^2)$ satisfying $F(z_k,w)=G(z_k,w)$ and $F(z,w_k)=G(z,w_k)$, for
$k=1,\dots, n$, with $\|\bar\bee F\|_{\BMOr/}\le 1$, as well as two
bounded functions on $\D^2$, $F_1$ analytic in $z$ and $F_2$ analytic
in $w$, satisfying $F_1(z_k,y)=G(z_k,y)$ and $F_2(x,w_k)=G(x,w_k)$,
for $k=1,\dots,n$, with $\|F_1\|_\infty\le1$ and $\|F_2\|_\infty\le1$,
whenever the two numerical $n\times n$ matrices
$$
((1-\|G_2\|_\infty^2G_1(z_j)\overline{G_1(z_k)})(1-z_j\overline{z_k})^{-1})_{j,k=1,\dots,n}
$$
and
$$
((1-\|G_1\|_\infty^2G_2(w_j)\overline{G_2(w_k)})(1-w_j\overline{w_k})^{-1})_{j,k=1,\dots,n}
$$
are positive definite.
\end{corollary}

\section{Hankel Operators of Finite Type and Versions of the Kronecker and
AAK Theorems}

In the one-dimensional case, once the relation between the bounded
Hankel operators and their symbols was established, it was important
to characterize the symbols of operators of finite rank.  The
characterization is given by the Kronecker theorem: A bounded Hankel
operator $\Gamma $ is of finite rank if and only if $\Gamma = \Gamma
_\varphi $ for $\varphi = \bar b h$, where $b$ is a finite Blaschke
product and $h\in H^\infty $, so $ \bar b h\in H^\infty + R_n$ (where
$R_n$ is the class of rational functions with $n$ poles in the disk).
Since the range of $\Gamma $ is finite-dimensional if and only if its
kernel has finite codimension, and since this kernel is a subspace of
$H^2(\T)$ invariant under the shift $S$, the Kronecker theorem can be
deduced from the Beurling theorem, asserting that a subspace ${\funnyT}
\subset H^2(\T)$ is invariant if and only if ${\funnyT} = \theta
H^2(\T)$, where $\theta $ is an inner function with $|\theta | \equiv
1$, and that an invariant subspace ${\cal \T}$ is of finite codimension
if and only if $\theta = b$, a finite Blaschke product.  The
$S^*$-invariant subspaces, called the model spaces, are of the form
$K_\theta = H^2(\T) \ominus \theta H^2(\T)$, and $K_\theta $ is
finite-dimensional if and only if $\theta = b.$

Recall that for an operator $T$ and for $n\in \N$, the singular
numbers of $T$ are defined as
\begin{equation}
s_n(T) := \inf  \{\|T - T_n\|: T_n\hbox{ of finite rank } \leq  n\},
\label{eq3.1}
\end{equation}
which is equivalent to
$$
s_n(T) = \inf  \{\|T \big| _E\|: E\hbox{ of codimension } \leq  n\}.
$$
Here $s_0 (T) = \|T\| \geq s_1 \geq \cdots\geq s_n\geq\cdots$, and $T$
is of finite rank if there is an $m\in \N $ such that
$s_n(T) = 0$ for $n > m$.

A theorem of Adamjan--Arov--Krein [AAK] asserts that for every Hankel
operator $\Gamma : H^2(\T)
\rightarrow  H^2(\T)^\bot $ we have, for $n\in \N $,
\begin{equation}
s_n(\Gamma ) = \inf  \{\|\Gamma  - \Gamma _n\|: \Gamma _n
\hbox{ Hankel and of finite rank $\leq  n\}$}.
\label{eq3.2}
\end{equation}
This, combined with the Kronecker theorem, gives, for all
$\varphi \in L^\infty (\T)$,
$$
s_n(\Gamma _\varphi ) = \dist(\varphi,\,H^\infty  + R_n).
$$

An equivalent form of \eqref{3.2} was given by S. Treil [T1] as
\begin{equation}
s_n(\Gamma ) = \inf \{\|\Gamma |_{\funnyT} \|: {\funnyT}
\hbox{  invariant under $S$ and codim }{\funnyT}  \leq  n\},
\label{eq3.2a}
\end{equation}
where, by Beurling's theorem, the subspace $\funnyT$
is of the form $bH^2(\T)$,
for $b$ a Blaschke product with $n$ factors.

Through an abstract version of the AAK theorem, in [CS1] it was shown that
the situation is radically different for (big) Hankel operators:

{\def\thelemma{B}\addtocounter{lemma}{-1}
\begin{theorem}
For every bounded $\Gamma : H^2(\T^2) \rightarrow
H^2(\T^2)^\bot $ and for all $n\in \N $, we have
$$
s_n(\Gamma ) \geq  1/\sqrt{2} \|\Gamma \|.
$$
\end{theorem}}

This theorem implies that all (big) Hankel operators of finite rank are zero,
and no satisfactory extension of the AAK theorem can  be expected in terms
of their singular numbers.  This is linked to the fact that, by a theorem of
Ahern and Clark [AhC], the subspaces of the form $bH^2(\T^2)$
for $b = b_1 \otimes  b_2$, with 
$b_1$ and $b_2$ finite one-dimensional Blaschke
products, are not of finite codimension in $H^2(\T^2)$. However, as shown in
[CS4], for these subspaces $bH^2(\T^2)$, it is still true that 
$$
\{f_{y_0}(x) =f(x,y_0):f(x,y) \in bH^2(\T^2), \,y_0 \in \T \hbox{ fixed}\}
$$
and
$$
\{f_{x_0}(y) =f(x_0,y):f(x,y) \in bH^2(\T^2), \,x_0 \in \T \hbox{fixed}\}
$$ 
are finite-codimensional subspaces of $H^2(\T)$, leading to a notion
of subspaces of finite bi-codi\-mension.

A decomposable subspace $V$ of $H^2(\T^2)$, with $V = V_1 \otimes
V_2$, where $V_1 \subset H^2(\T)$ and $V_2 \subset H^2(\T)$, is called
of {\it finite bi-codimension} $(m,n)$ if and only if $\codim V_1=m$
and $\codim V_2=n$.
The orthogonal complement of such $V$ is
$$
V^\perp=H^2(\T^2)\ominus V=V_1^\perp\otimes H^2(\T)+H^2(\T)\otimes
V_2^\perp,
$$
with $V_k^\perp=H^2(\T)\ominus V_k$ for $k=1,2$.  Since orthogonal
complements of this form will appear again in Section 5 and in other
contexts, we give them a name.  Given a subspace $\cal L\subset
L^2(\T)$ and two positive integers $m$ and $n$, a subspace $W\subset
L^2(\T^2)$ is said to be of {\it bi-finite type\/} $(\cal L;m,n)$ if
and only if there exist two subspaces $W_1$ and $W_2$ of $H^2(\T)$,
with $\dim W_1=m$ and
$\dim W_2=n$, such that $W=W_1\otimes\cal L+\cal L\otimes W_2$.  With
this notation, the orthogonal complement $V^\perp$ of a subspace
$V\subset H^2(\T^2)$ of finite bi-codimension $(m,n)$ is a subspace of
bi-finite type $(H^2(\T);m,n)$.  

For these notions we have the following analogue to the Beurling theorem for
invariant subspaces of finite codimension in the disk:

{\def\thelemma{C}\addtocounter{lemma}{-1}
\begin{theorem}{\normalfont[CS4]}
For a subspace ${\funnyT}  \subset  H^2(\T^2)$,
invariant under both shifts, $S_1$ and $S_2$, of $H^2(\T^2)$,
the following conditions are equivalent:
\begin{enumerate}
\item[\normalfont(a)] 
${\funnyT}$  is of finite bi-codimension $(m,n)$, that is,
$\funnyT=V_1\otimes V_2$;
\item[\normalfont(b)] 
$\funnyT^\perp=H^2(\T^2)\ominus\funnyT$ is of bi-finite type
$(H^2(\T);m,n)$, that is, $\funnyT^\perp=W_1\otimes H^2(\T)+H^2(T)\otimes
W_2$ with $\dim W_1=m$ and $\dim W_2=n$.
\item[\normalfont(c)] 
${\funnyT} = bH^2(\T^2) = b_1H^2(\T) \otimes b_2H^2(\T)$, for $b = b_1
\otimes b_2$, that is, $b(x,y) = b_1(x)b_2(y)$, where $b_1$ and $b_2$
are one-dimensional Blaschke products with $m$ and $n$ factors,
respectively.
\end{enumerate}
Furthermore, $W_k=V_k^\perp=H^2(\T)\ominus b_kH^2(\T)$ for $k=1,2$.
\end{theorem}}

A Hankel operator $\Gamma : H^2(\T^2) \rightarrow  H^2(\T^2)^\bot $ is
called of {\it finite type $(m,n) \in \N ^2$} if the kernel of
$\Gamma $ is of finite bi-codimension $(m,n).$

\begin{theorem}{\bf (Kronecker-type characterization of Hankel operators of
finite type)}
\label{3.1}
Let $\Gamma : H^2(\T^2) \rightarrow  H^2(\T^2)^\bot $ be a bounded Hankel
operator.  The following conditions are equivalent:

\begin{enumerate}
\item[\normalfont(a)] 
$\Gamma $ is of finite type $(m,n)$.
\item[\normalfont(b)] 
$\Gamma = \Gamma _\phi $ for $\phi \in \BMOr/ \cap H^{2\bot }$, and
$\phi = \bar b h$, where $b = b_1 \otimes b_2$ for $b_1$ and $b_2$
one-dimensional Blaschke products with $m$ and $n$ factors,
respectively, and $h\in H^2(\T^2)$.  Moreover, $h = bh_x + \varphi _1 =
bh_y + \varphi _2$ for $\varphi _1,\varphi _2 \in L^\infty (\T^2)$,
$h_x \in H^2_x(\T^2)$, $h_y \in H^2_y(\T^2)$, so that
$$
\phi  = h_x + \bar b \varphi _1 = h_y + \bar b \varphi _2.
$$
\end{enumerate}
\end{theorem}

\begin{proof}  
The kernel $K$ of $\Gamma $ is invariant under both shifts $S_1$ and
$S_2$, since $\Gamma $ is Hankel and, for $f\in K$, $\<\Gamma S_kf,g\>
= \<\Gamma f$, $S^{-1}_kg\> = 0$, for $k = 1,2$.  Then, by Theorem C,
there is $b = b_1 \otimes b_2$ such that $K = bH^2(\T^2)$.  For $\phi
\in \BMOr/$ the $H^{2\bot }$ symbol of $\Gamma $, we have $\<\Gamma
f,g\> = \int f\bar g\phi $, and thus $\int bf\bar g\phi = 0$ for all
$f\in H^2$ and $g\in H^{2\bot }$, that is, $b\phi = h\in H^2$ and
$\phi = \bar bh$, and the converse holds.  Moreover, since $\bar bh =
\phi \in \BMOr/$, the function $h$ must satisfy $h = bh_x + \varphi _1
= bh_y + \varphi _2$ for $\varphi _1,\varphi _2 \in L^\infty (\T^2)$,
$h_x \in H^2_x$, $h_y \in H^2_y$, and the conclusion follows.
\end{proof}

Theorem C suggests, in order to develop a version of the AAK theorem, to
replace the ordinary singular numbers of Hankel operators by some
$\sigma $-numbers defined in analogy with \eqref{3.2a}. For $\Gamma : H^2(\T^2)
\rightarrow  H^2(\T^2)^\perp $ a Hankel operator and $(m,n)
\in \N ^2$, let
\begin{equation}
\sigma _{mn}(\Gamma ) := \inf \{\|\Gamma |_{\funnyT} \|: {\funnyT} \subset
H^2(\T^2)\vtop{\hbox{ invariant under $S_1$ and $S_2$}\hbox{ and
of finite bi-codimension $(m,n)\}$.}}
\label{eq3.4}
\end{equation}
Equivalently, for $\phi \in  \BMOr/$,
{
\thickmuskip=.7 \thickmuskip
\medmuskip=.7 \medmuskip
\begin{equation}
\sigma _{mn}(\Gamma _\phi ) = \inf \{\|\Gamma _{b\phi }\|: b = b_1 \otimes
b_2, \vtop{\hbox{with $b_1$ and $b_2$ one-dimensional Blaschke products}
\hbox{having at most $m$ and $n$ factors, respectively\}}},
\label{eq3.4a}
\end{equation}}%
since, by Theorem C, the subspaces ${\funnyT}$  in \eqref{3.4} can be written as
${\funnyT} = bH^2(\T^2)$ with $\|bh\|_2 = \|h\|_2$, so that
$$
\|\Gamma _\phi |_{\funnyT} \| = \sup_h
 {\|\Gamma _\phi  bh\|_2\over \|bh\|_2} = \sup_h
 {\|\Gamma _{b\phi }h\|_2\over \|h\|_2} = \|\Gamma _{b\phi }\|.
$$

It is easy to check that the infima are attained in \eqref{3.4} and \eqref{3.4a}.  Clearly, we
have again
$$
\eqalign{
\sigma _{00}(\Gamma ) &= \|\Gamma \|,\cr
\sigma _{mn}(\Gamma ) &\geq \sigma _{(m+1)n}(\Gamma ), \cr
\sigma _{mn}(\Gamma ) &\geq \sigma _{m(n+1)}(\Gamma ),\cr
\sigma _{mn}(\Gamma ) &\leq  \|\Gamma \|}
$$
for all $m,n\in \N $.

\begin{corollary}
For every $(m,n) \in \N ^2$ there exists
a non-zero Hankel operator $$\Gamma : H^2(\T^2)
\rightarrow  H^2(\T^2)^\bot$$ of finite type $(m,n)$, such that
$$
\sigma _{pq}(\Gamma ) = 0\hbox{ for }p > m,\quad q > n.
$$
\end{corollary}

\begin{proof}  
Take $h(x,y) = b_1(x)h_1(y) + b_2(y)h_2(y)$ for $b_1$ and $b_2$
one-dimensional Blaschke products with at most $m$ and $n$ factors,
respectively, and $h_1,h_2 \in H^\infty (\T)$.  Further take $\phi =
\bar bh$, for $b = b_1 \otimes b_2$.  Then $h\in H^\infty (\T^2)$, and
$$
\phi  = \overline{ b_1(x)} h_2(x) + \overline{ b_2(y)} h_1(y)
\in  \BMOr/ \cap  H^2(\T^2)^\bot,
$$
since 
$$
P_{-x}\phi  = P_{-x}\overline{b_1(x)}
h_2(x), \qquad
P_{-y}\phi  = P_{-y}\overline{b_2(x)}
h_1(y),
$$ 
with $ \bar b_1 h_2$ and $ \bar b_2h_1$ in $L^\infty (\T^2)$, and
$P_xP_y\phi = 0$ for the right choice of $h_1$ and $h_2$.  By
Theorem~\eqref{3.1}, the Hankel operator $\Gamma = \Gamma _\phi $
satisfies the conclusion. 
\end{proof}

From Theorem B follows that there are no nonzero compact big Hankel
operators, that is,
Hankel operators whose sequence of singular numbers tend to zero.
Since this corollary says that there are big Hankel operators $\Gamma  \neq  0$
with $\sigma _{mn}(\Gamma ) \rightarrow  0$ as $m,n \rightarrow  \infty $, it
is interesting to study the class of such operators, and this will be done elsewhere.

For bounded Hankel operators in the one-dimensional case the AAK theorem
asserts that $s_n(\Gamma _\varphi ) = \dist_{L^\infty }(\varphi, H^\infty  +
R_n)$.

This precise statement does not hold for all bounded (big) Hankel operators
$\Gamma  = \Gamma _{\phi '}$ given by a symbol $\phi \in  \BMOr/$, but we
still have a substitute by replacing the distance
{\thickmuskip=.6 \thickmuskip
\medmuskip=.6 \medmuskip
\hfuzz=2pt
$$
\eqalign{\dist (\phi,\, H^\infty  + R_n) &= \inf \{\|\phi- b^{-1}h\|_\infty :
h\in H^\infty  = L^\infty  \cap  H^2,\cr&\qquad  b = b_1 \otimes 
\cdots\otimes b_d,\ \vtop{\hbox{with $b_k$ a one-dimensional Blaschke product}
\hbox{of at most $n_k$ factors, for $k = 1,\dots, d$\}}}\cr
&= \inf \{\|b\phi  - h\|_\infty : h\in H^\infty  = L^\infty  \cap
H^2,b\}}
$$ 
by
$$
\eqalign{\delta (\phi,\, \BMOAr/  + R_n) &:= \inf \{\|b\phi  - h\|_{\BMOr/}:
h\in \hbox{ \BMOAr/ } = \BMOr/ \cap  H^2,\cr&\qquad b = b_1 \otimes 
\cdots\otimes b_d, \vtop{\hbox{with $b_k$ a one-dimensional Blaschke product}
\hbox{of at most $n_k$ factors, for $k = 1,\dots, d$\}.}}}
$$}%
Observe that $\BMOAr/=\BMOr/\cap H^2=\BMO/\cap H^2=\BMOA/$.

\begin{theorem}  
\label{3.3}
For every $\phi \in  \BMOr/(\T^d)$ and $n\in \N ^d$, with $d > 1$, we have
$$
1/\sqrt{d}\sigma _n(\Gamma _\phi ) \leq  \delta (\phi,\, \BMOAr/  + R_n)
\leq  \sigma _n(\Gamma _\phi ),
$$
where $\Gamma _\phi $ is the Hankel operator with symbol $\phi$.
\end{theorem}

\begin{proof}
By \eqref{3.4a}, for every $\epsilon > 0$, there are $b_1$ and $b_2$,
and $b = b_1 \otimes b_2$, such that
$$
\sigma _{mn}(\Gamma _\phi ) \leq  \|\Gamma _{b\phi }\| \leq
\sigma _{mn}(\Gamma _\phi ) + \epsilon .
$$
The operator $\Gamma _{b\phi }: H^2 \rightarrow H^{2\bot }$ is also
Hankel, and, by Theorem A, $\Gamma _{b\phi } = \Gamma _\psi $ for some
$\psi \in \BMOr/$, with $\psi = b\phi - h$, $h\in H^2$, and $\|\Gamma
_{b\phi }\| \leq \sqrt{2}\|\psi \|_{\BMOr/}$, $\|\psi \|_{\BMOr/} \leq
\sigma _{mn}(\Gamma _\phi ) + \epsilon $.  Hence, $1/\sqrt{2} \sigma
_{mn}(\Gamma _\phi ) \leq \|b\phi - h\|_{\BMOr/} \leq \sigma
_{mn}(\Gamma _\phi ) + \epsilon $, for all $\epsilon > 0$, which is the
conclusion.
\end{proof}

\section{Carleson Measures, Model Subspaces of Finite Type, and \BMOr/ Symbols}

In the one-dimensional case, there is a close relation linking Hankel
operators, \BMO/ functions and Carleson measures.  Carleson measures
in the disk are those positive measures $\mu $ satisfying the Carleson
imbedding condition
\begin{equation}
\label{eq4.1}
\int  _\D|f(z)|^2 \,d\mu (z) \leq  C^2 \int  _\T|f(t)|^2\,dt,\quad
\hbox{for all }f\in H^2(\T),
\end{equation}
where $f(z)$ stands for the analytic extension of $f$ to $\D$.

Carleson characterized those measures as satisfying the {\it tent
condition\/} for intervals, that is, $\mu(S(I))\le C|I|$ for every
interval $I$, where $S(I)$ is a tent in $\D$ with base $I$.  Moreover,
the $H^2$-imbedding condition \eqref{4.1} is equivalent to the
$H^p$-imbedding condition being valid for all $p$ such that $1\le
p<\infty$.

Following Nikolskii and Treil [Ni], [T2], condition \eqref{4.1} can be
expressed in terms of projectors on one-dimensional model subspaces
$$
K_z = K_{b_z} = H^2 \ominus  b_zH^2
$$
defined by single-factor Blaschke products
$$
\label{eq4.2a}
b_z(\zeta ) = {|z|\over z} {z - \zeta \over 1 - \zeta \bar z}.
$$

It is well known that such a subspace $K_z$ is spanned by the normalized
function
\begin{equation}
\label{eq4.2b}
\phi _z(\xi ) = {(1 - |z|^2)^{1/2}\over 1 - \bar z\xi }\quad\hbox{for
$z\in \D$ and $\xi \in \T$},
\end{equation}
which has the reproducing property
\begin{equation}
\label{eq4.2c}
\<f,\phi _z\> = (1 - |z|^2)^{1/2} f(z)\quad\hbox{for all $f\in H^2$}.
\end{equation}
Thus, for $P_z: H^2(\T) \rightarrow  K_z$ the orthogonal projector, the
identity
\begin{equation}
\label{eq4.3}
\|P_zf\|^2_2 = (1 - |z|^2) |f(z)|^2
\end{equation}
holds for all $f\in H^2(\T)$.  The Carleson imbedding condition \eqref{4.1}
can, therefore, be rewritten as
\begin{equation}
\label{eq4.1a}
\int  _\D\|P_zf\|^2_2 \,d\nu (z) \leq  C^2 \|f\|^2_2,\quad\hbox{for
all $f\in H^2(\T)$},
\end{equation}
with $d\nu (z) = (1 - |z|^2)^{-1} \,d\mu (z).$

Moreover, for $P: L^2 \rightarrow  H^2$ the analytic projector, we have
\begin{equation}
\label{eq4.3a}
P_zf = b_z(I - P)\bar b_zf = b_z\Gamma_{\bar b _z}f
\end{equation}
and
\begin{equation}
\label{eq4.3b}
\|P_zf\|^2_2 = \|\Gamma_{\bar b _z}f\|^2_2.
\end{equation}

From \eqref{4.3b} it can be deduced (see the development leading to
\eqref{4.9} below) that $\mu \geq 0$ is Carleson if and only if a
canonically associated (vector-valued) Hankel operator $\Gamma $ is
bounded, and (through the Nehari theorem) if and only if its
antianalytic symbol $\Gamma 1\in \BMO/$.

In Section 1 we observed that the different definitions of $\BMO/$,
which coincide for $d=1$, give rise to different classes in $\T^d$,
for $d>1$.  In fact, Chang and Fefferman defined product \BMO/ to
circumvent Carleson's counterexample showing that the class of
measures in $\D^2$ characterized by the tent condition on rectangles
$R=I\times J$ does not necessarily satisfy the $H^1$-imbedding
condition.  Enlarging the class of tents in their definition of
product Carleson measures, Chang and Fefferman proved that a function
is in product $\BMO/$ if and only if a canonically associated measure
is product Carleson.  Here we adopt in $\D^d$ the Nikolskii
formulation \eqref{4.1a} and show that a measure is
Carleson--Nikolskii if and only if a canonically associated function
is $\BMOr/$.

In $\T^2$ the one-dimensional subspace $K_{b_z}$, where $b_z$ is a
Blaschke factor, is replaced by $K_{b_z\otimes b_\zeta } =: K_{z\zeta
}$, where $b_z$ and $b_\zeta $ are one-variable Blaschke factors, and
$K_{z\zeta }$ is not a one-dimensional subspace of $H^2(\T^2)$.  But
now, according to Theorem C in Section 4, $K_{z\zeta }$ is of
bi-finite type $(H^2(\T);1,1)$, and its elements are of the form $A(y)\phi
_z(x) + B(x)\phi _\zeta (y)$, for $A(y)$
and $B(x)$ varying in $H^2(\T)$.  For the orthogonal projector $P_{z\zeta }$
from $H^2(\T^2)$ onto $K_{z\zeta }$ we have the following result.

\begin{lemma}[Lemma on the projection] 
\label{4.1}
If $P_{z\zeta }: H^2(\T^2)
\rightarrow  K_{z\zeta }$ is the orthogonal projector, we have, for all
$f\in H^2(\T^2)$,
\begin{equation}
\label{eq4.4}
P_{z\zeta }f(x,y) = c_zf(z,y)\phi _z(x) + c_\zeta f(x,\zeta )\phi _\zeta (y) -
c_zc_\zeta f(z,\zeta )\phi _z(x)\phi _\zeta (y),
\end{equation}
where $c_\nu  = (1 - |\nu |^2)^{1/2}$, $\nu \in \D$.
Equivalently,
$$
P_{z\zeta }f(x,y) = \<f,\phi _z\>_{L^2(\T_x)}\phi _z(x) +
\<f,\phi _\zeta \>_{L^2(\T_y)}\phi _\zeta (y)
- \<f,\phi_z  \otimes  \phi _\zeta \>_{L^2(\T^2)}\phi _z(x)\phi _\zeta (y).
$$
\end{lemma}

\begin{proof}  
Denoting the right-hand side of \eqref{4.4} by $g(x,y)$, for 
$g\in K_{z\zeta }$, it remains to check that, for arbitrary
$A,B\in H^2(\T)$, we have
$$
\eqalign{
\<g(x,y),\,A(y)\phi _z(x)\> &= \<f(x,y), \,A(y)\phi _z(x)\>,\cr
\<g(x,y),\,B(x)\phi _\zeta (y)\> &= \< f(x,y), \,B(x)\phi _\zeta (y)\>.}
$$
Since, by \eqref{4.2c}, for every $F\in H^2(\T^2)$ we have
\begin{equation}
\label{eq4.5}
\<F(x,y), \,\phi _z(x)\> = c_zF(z,y)\quad\hbox{and}\quad
\<F(x,y), \,\phi _\zeta (y)\> = c_\zeta F(x,\zeta ),
\end{equation}
and, by \eqref{4.2b}, $\int   |\phi _z(x)|^2\,dx = \int   |\phi _\zeta
(y)|^2\,dy = 1$, we obtain, as desired,
$$
\eqalign{
\<g(x,y), \,A(y)\phi _z(x)\> &= c_z\int   f(z,y)
\overline{A(y)}
\,dy + c_\zeta \iint    f(x,\zeta )\phi _\zeta (y)
\overline{A(y)} \overline{\phi_z(x)} \,dx\,dy
\cr&\qquad
- c_zc_\zeta f(z,\zeta ) \iint    \phi _z(x)\phi _\zeta (y)
\overline{A(y)} \overline{\phi_z(x)} \,dx\,dy\cr
&= c_z\int   f(z,y) \overline{A(y)}\, dy + c_\zeta c_zf(z,\zeta )c_\zeta
\overline{A(\zeta)} - c_zc_\zeta f(z,\zeta )c_\zeta
\overline{A(\zeta)} \cr&= \<f(x,y), \,A(y)\phi_z(x)\>,}
$$
and similarly for the other term. 
\end{proof}

\begin{lemma}
\label{4.2}
For $P_{z\zeta }: H^2(\T^2) \rightarrow  K_{z\zeta }$, where
$(z,\zeta ) \in \D^2$, and $f\in H^2(\T^2)$, we have
$$
\eqalign{
\|P_{z\zeta }f\|^2_2(1 - |z|^2)^{-1} (1 - |\zeta |^2)^{-1} &= (1 -
|\zeta |^2)^{-1} \int  _\T|f(z,y)|^2\,dy\cr&\qquad
+ (1 - |z|^2)^{-1} \int  _\T|f(x,\zeta )|^2\,dx - |f(z,\zeta )|^2,}
$$
where $f(z,y)$, $f(x,\zeta )$ and $f(z,\zeta )$ are the analytic
extensions of $f$ to $z\in \D$, $\zeta \in \D$,
and $(z,\zeta ) \in \D^2$, respectively.
\end{lemma}

\begin{proof}  
From the expression of $P_{z\zeta }f$ given in \eqref{4.4} it follows, using
\eqref{4.2b}, \eqref{4.2c}, \eqref{4.3} and \eqref{4.5}, that
$$
\displaylines{ |P_{z\zeta }f|^2 = c^2_z|f(z,y)|^2 |\phi _z(x)|^2 +
c^2_\zeta |f(x,\zeta )|^2 |\phi _\zeta (y)|^2 + c^2_zc^2_\zeta
|f(z,\zeta )|^2 |\phi _z(x)|^2 |\phi _\zeta (y)|^2 \hfill\cr\hfill{} +
2\Re\Bigl(c_zc_\zeta f(z,y) \overline{f(x,\zeta)} \phi _z(x)
\overline{\phi_\zeta(y)} + c^2_zc_\zeta \overline{f(z,\zeta)}
f(z,y)|\phi _z(x)|^2 \phi _\zeta (y)\hfill\cr\hfill - c_zc^2_\zeta
f(z,\zeta ) \overline{f(x,\zeta)} \phi _z(x)|
\overline{\phi_\zeta(y)}|^2\Bigr).}
$$
Then,
$$
\thickmuskip=.5\thickmuskip
\medmuskip=.5\medmuskip
\eqalign{
\|P_{z\zeta }f\|^2_2 &= \iint    |P_{z\zeta }f(x,y)|^2\,dx\,dy \cr&=
c^2_z\int   |f(z,y)|^2\,dy + c^2_\zeta \int   |f(x,\zeta )|^2\,dx
+ c^2_zc^2_\zeta |f(z,\zeta )|^2 - 2\Re(c^2_zc^2_\zeta |f(z,\zeta)|^2),}
$$
so
$$
(c^2_zc^2_\zeta )^{-1} \|P_{z\zeta }f\|^2_2 = c^{-2}_\zeta  \int
|f(z,y)|^2\,dy + c^{-2}_z \int   |f(x,\zeta )|^2\,dx - |f(z,\zeta )|^2,
$$
which is the conclusion. 
\end{proof}

Following Nikolskii's approach, we say that a measure $\mu  \geq  0$ defined
in $\D^2$ is {\it Carleson--Nikolskii} if
$$
d\nu (z,\zeta ) = (1 - |z|^2)^{-1} (1 - |\zeta |^2)^{-1} d\mu (z,\zeta )
$$
satisfies
$$
\iint   _{\D^2}\|P_{z\zeta }f\|^2_2 \,d\nu (z,\zeta ) \leq  C^2
\|f\|^2_2\quad  \hbox{for all $f\in H^2(\T^2)$}.
$$
Formula \eqref{4.3a} is still valid in $H^2(\T^2)$, that is,
for all $f\in H^2(\T^2)$ we have
$$
P_{z\zeta }f = (b_z \otimes  b_\zeta ) (I - P)(\bar b_z \otimes
\bar b _\zeta )f = (b_z \otimes  b_\zeta )\Gamma _{
\bar b _z\otimes \bar b _\zeta}f.
$$
Thus, again we have
$$
\|P_{z\zeta }f\|^2_2 = \|\Gamma _{\bar b _z \otimes 
\bar b _\zeta }f\|^2_2.
$$
Therefore, a measure $\mu $ is Carleson--Nikolskii if and only if
$$
\iint   _{\D^2}\|\Gamma_{ \bar b _z \otimes 
\bar b _\zeta} f\|^2_2 \,d\nu  \leq  C^2 \|f\|^2_2\quad
\hbox{for all }f\in H^2(\T^2).
$$

Let $\Gamma : H^2(\T^2) \rightarrow  L^2(\D^2,\nu ; H^2(\T^2)^\perp )$ be
the operator assigning to each $f\in H^2(\T^2)$ the function
$$
(z,\zeta ) \mapsto \Gamma_{
\bar b _z\otimes \bar b _\zeta} f\in H^2(\T^2)^\perp,
$$
so that $\mu $ is of Carleson type if and only if $\Gamma $ is bounded in
$L^2(\D^2,\nu )$, with $\|\Gamma \| \leq  C$.  By Fubini's theorem, the
space $L^2(\D^2,\nu ; H^2(\T^2)^\perp )$ of square integrable functions
in the bidisk, with values in $H^2(\T^2)^\perp $, is isometrically isomorphic to
the space $H^2(\T^2)^\perp  (L^2(\D^2,\nu ))$ of antianalytic functions
with values in $L^2(\D^2,\nu )$.  Under this isomorphism the operator
$\Gamma $ corresponds to the operator
\begin{equation}
\vec\Gamma 
: H^2(\T^2) \rightarrow  H^2(\T^2)^\perp  (L^2(\D^2,\nu )).
\label{eq4.9}
\end{equation}
This $\vec\Gamma$ is a (vector-valued) big Hankel operator, since, for
$k = 1,2$, we have
$$
\vec\Gamma S_kf = F(x,y; z,\zeta ) = \Gamma_{\bar b _z \otimes \bar b
_\zeta} (S_kf)(x,y) = (I - P)S_k \Gamma_{\bar b _z\otimes \bar b
_\zeta}f = (I - P)S_k \vec\Gamma f.
$$
The operator $\vec\Gamma$ is called {\it the Hankel operator
canonically associated to $\mu $}.  Theorem A (which can be used since
its proof through abstract liftings extends to Hankel operators from
the scalar spaces $H^2(\T^2)$ to a vector-valued $H^2(\T^2)^\perp ({\cal
H} )$, where ${\cal H}$ is a Hilbert space) applied to $\vec\Gamma$
yields:

\begin{theorem}  
\label{4.3} 
A measure $\mu \geq 0$ in $\D^2$ is of Carleson type, with constant
$C$, if and only if the canonically associated operator $ \vec\Gamma $
is bounded with norm $\| \vec\Gamma \| = C$, and if and only if
$\vec\Gamma 1\in \BMOr/(\T^2; L^2(\D^2,\nu ))$, with norm $\cong C$.
\end{theorem}  

Thus the connection between measures satisfying the Carleson imbedding
condition, Hankel operators and $\BMO/(\T)$, is recovered in $\T^2$ in terms of
\BMOr/.

\section{Estimates for the Norm of the Hankel Operators of Finite Type} 

Let us recall some basic properties of the finite-dimensional model
subspaces $K_b \subset H^2(\T)$, where $b$ is a finite Blaschke
product, which include as a special case the properties of the $K_z$
considered in Section 5.  For a finite Blaschke product $b$ with
simple zeros in $\D$, we again denote by $P_b$ the orthogonal
projector from $H^2(\T)$ onto $K_b = H^2(\T)\ominus bH^2(\T)$, and
define the model operator $T_b: K_b \rightarrow K_b$ by $T_b :=
P_bS|K_b$, so that $T^*_b = S^*|K_b$.  Similarly, for each $G\in
H^\infty (\T)$, $G(T_b)$ is defined by $G(T_b)f := P_bGf.$

If $\phi _z$ is given by \eqref{4.2b}, then, for each
$z\in \D$, $\phi _z$ is an eigenfunction of $S^*$,
and if $z_1,\dots, z_m \in \D$ are the
zeros of $b$, then $\{\phi _{z_1},\dots, \phi _{z_m}\}$ is a basis of $K_b$
composed of eigenfunctions of $T^*_b$.  Similarly, $K_b$ has a basis
$\{\psi _{z_1},\dots, \psi _{z_m}\}$, of eigenfunctions of $T_b$, where
\begin{equation}
\label{eq6.1}
\psi _z(\xi ) = b(\xi  - z)^{-1},\quad T_b\psi_z=z\psi_z.
\end{equation}
Thus $T_b$ and $T^*_b$ are {\it multiplier operators}, that is,
they are given by
diagonal finite matrices in the corresponding bases, so that, for each
$G\in H^\infty (\T)$, the condition
$$
\|G(T_b)\| \leq  1
$$
is equivalent to the positive definiteness of the associated Pick matrix
$$
((1 - G(z_j)\overline{G(z_k)}) (1 - z_j\overline{z_k})^{-1})_{j,k=1,\dots,m}.
$$

The Kronecker theorem characterizes the symbols of the Hankel
operators $\Gamma:H^2(\T)\to H^2(\T)^\perp$ of finite rank $n$ as
those of the form $\bar b G$, for $b$ a Blaschke product with $n$
factors and $G\in H^\infty(\T)$.  Since the projector $P_b$ is
related to the analytic projector $P:L^2\to H^2$ by
\begin{equation}
P_bf=b(I-P)\bar b f=b\Gamma_{\bar b}f,
\label{eq5.4}
\end{equation}
we derive the identities
$$
|\Gamma_{\bar bG}f|=|\Gamma_{\bar b}Gf|=|P_bGf|=|G(T_bf)|,
$$
and thus
$$
\|\Gamma_{\bar bG}\|=\|G(T_b)\|.
$$
This means that, in the circle, {\it the norm of a Hankel operator of
finite rank is equal to the norm of an associated multiplier
operator\/} acting in finite-dimensional $K_b$, which in turn is
determined by a finite Pick matrix. 

The same result holds for Hankel operators of finite type in the torus
(see Theorem~\ref{5.3} below), but the association with the multiplier
operators acting in $K_b$ is not so simple.  This is due to the fact
that here $K_b$, for $b$ the tensor product of $d$ Blaschke products,
is not finite-dimensional but of multiple-finite type.  As before, we
present here the case $d=2$.

In $\T^2$, if we restrict ourselves to the case when $b_1$ and $b_2$ have the
same number of zeros, at $z_1,\dots, z_n$ and $w_1,\dots, w_n$, respectively,
and when $G = G_1 \otimes  G_2$, with $G_1,G_2 \in H^\infty (\T)$, we have the
following equivalences, in terms of
\begin{equation}
\label{eq5.9}
K_{b_1b_2} = H^2(\T^2) \ominus  (b_1 \otimes  b_2) H^2(\T^2)
\end{equation}
and
\begin{equation}
\label{eq5.9a}
K'_{b_1,b_2} =\overline{H^2(\T^2)^\perp}
 \ominus  [(b_1 \otimes  b_2) H^2(\T^2) \oplus  b_1(H^2_x \cap  H^2_{-y})
\oplus  b_2(H^2_{-x} \cap  H^2_y)]
\end{equation}
subspaces of finite type $(H^2(\T); n,n)$ and $(L^2(\T); n,n)$, respectively,
as follows from Theorem C and from Theorem 2 in [CS4] (see Section 4).

\begin{proposition}  
\label{5.2} Given $\phi = (\bar b_1 \otimes \bar b_2) (G_1 \otimes
G_2)$, for $b_1,b_2$, $G_1$ and $G_2$ as above, and $K_{b_1b_2}$ and
$K'_{b_1b_2}$ defined by \eqref{5.9} and \eqref{5.9a}, the following
conditions are equivalent:
\begin{enumerate}
\item[\normalfont(a)]
$\|\Gamma _\phi \| \leq  1$, that is, $\|\Gamma _\phi f\|_{L^2} \leq
\|f\|_{L^2}$, for all $f\in H^2(\T^2)$.
\item[\normalfont(b)]
$\|\Gamma_\phi|K_{b_1b_2}\|\le 1$, that is,
$\|\Gamma _\phi e\|_{L^2} \leq  \|e\|_{L^2}$,
for all $e\in K_{b_1b_2}$.
\item[\normalfont(c)]
For all $e\in K_{b_1b_2}$ and $e'\in K'_{b_1b_2}$, the inequality
$$
\biggl|\iint    e\bar e  \phi \,dx\,dy\biggr| 
\leq  \|e\|_{L^2}\|e'\|_{L^2}
$$
holds.
\end{enumerate}
\end{proposition}

\begin{proof}
(a)$\iff$(b).  Every $f\in H^2(\T^2)$ can be written as an orthogonal sum
$$
f(x,y) = b_1(x)b_2(y)h(x,y) + e(x,y),\quad\hbox{ with $ h\in H^2(\T^2)$
and $e\in K_{b_1b_2}$},
$$
and
$$
\eqalign{ \Gamma _\phi (b_1 \otimes b_2)h &= (I -
P)(b_1(x)b_2(y)h(x,y) \overline{b_1(x)} \overline{b_2(y)}
G_1(x)G_2(y))\cr& = (I - P) (h(x,y) G_1(x)G_2(y)) = 0,}
$$
since $hG_1G_2 \in H^2(\T^2)$.  Thus, $\Gamma _\phi f = \Gamma _\phi
e$, $\forall f\in H^2(\T^2)$ and $\|\Gamma _\phi f\|_2 = \|\Gamma _\phi
e\|_2 \leq \|e\|_2 \leq \|f\|_2$, so (b) implies (a).  The converse
follows from $K_{b_1b_2} \subset H^2(\T^2)$.

\smallskip
\noindent
(a)$\iff$(c).  
Similar proof, observing the equivalence of (a) with
$$
\biggl|\iint    f(x,y)\overline{g(x,y)}
\phi (x,y)\,dx\,dy\biggr| \leq  \|f\|_2 \|g\|_2
$$
for all $f\in H^2(\T^2)$ and $g\in\overline{H^2(\T^2)^\perp}
$, and writing in terms of the decomposition of $H^2(\T^2)^\perp $ in the
direct sum of $K'_{b_1b_2}$ and its orthogonal complement. 
\end{proof}

Proposition~\ref{5.2} says that $\|\Gamma _\phi \| = \|\Gamma _\phi
|K_{b_1b_2}\|$, and we will prove that $\|\Gamma _\phi |K_{b_1b_2}\|$
coincides with the norm of a multiplier operator in $K_{b_1b_2}$ (see
Theorem~\ref{5.3} below), thus generalizing the one-dimensional
results.

The systems of eigenfunctions $\{\psi _{z_1},\dots, \psi _{z_n}\}$
and $\{\psi _{w_1},\dots, \psi _{w_n}\}$, where
$\{z_1,\allowbreak\dots,\allowbreak z_n\}$ and $\{w_1,\allowbreak\dots,\allowbreak w_n\}$ are the zeros of $b_1$ and $b_2$,
are bases for $K_{b_1}$ and $K_{b_2}$, respectively.  Through Theorem
C of Section 4, this
allows to write the elements $e\in K_{b_1b_2}$ as
\begin{equation}
\label{eq5.12}
e(x,y) = \sum ^n_{i=1} A_i(y) \psi _{z_i}(x) + \sum ^n_{j=1}
B_j(x)\psi _{w_j}(y)
\end{equation}
where, for $i,j = 1,\dots, n$, we have $A_i$, $B_j \in H^2(\T)$.  In
what follows we write, for simplicity,
$$
\psi _{z_i}(x) = \xi _i(x)\hbox{ and }\psi _{w_j}(y) = \eta _j(y)\quad
\hbox{for $i,j = 1,\dots,n$},
$$
and remark that the $\xi _i's$ are eigenfunctions of the model operator
$T_{b_1}$, for which,
by \eqref{6.1},
$T_{b_1}\xi _i = z_i\xi _i$, and the $\eta _j's$
are eigenfunctions of $T_{b_2}$, with
$T_{b_2}\eta _j = w_j\eta _j$.  For their part,
each $A_i(y)$, $B_j(x)$ can be written as
\begin{equation}
\label{eq5.14}
\eqalign{
A_i(y) &= b_2(y)h''_i(y) + \sum _k c_{ik}\eta _k(y),\cr
B_j(x) &= b_1(x)h'_j(x) + \sum _\ell  d_{j\ell } \xi e(x)\cr}
\end{equation}
with $h''_i, h'_j \in H^2(\T)$, for $i,j = 1,\dots, n$.  Moreover
$$
\Gamma _\phi e = (I - P)\phi e = (P_{-y} + P_{-x}P_y)\phi e = (P_{-x} +
P_xP_{-y})\phi e,
$$
and, by \eqref{5.4} and the definitions of $\T_b$ and $G(\T_b)$, we get
$$
\eqalign{
P_{-x}(\bar b_1G_1\xi _i)(x) &= \bar b_1T_{b_1}G_1\xi _i(x) =
\bar b_1G_1(z_i)\xi _i(x),\cr
P_{-y}(\bar b_2G_2\eta _j)(y) &= \bar b_2G_2(w_j)\eta _j(y).}
$$
Since $P_x\bar b_1G_1\xi _i = \bar b_1G_1\xi _i -
P_{-x}\bar b_1G_1\xi _i$, we have
$$
\eqalign{P_x\overline{b_1(x)} G_1(x)\xi _i(x) &= \overline{b_1(x)}
(G_1(x) - G_1(z_i))\xi _i(x),\cr
P_y\overline{b_2(y)}G_2(y)\eta _j(y) &= \overline{b_2(y)}
(G_2(y) - G_2(w_j))\eta _j(y).}
$$

Now  for every $e\in K_{b_1b_2}$, we can write $\Gamma _\phi e$ in
terms of functions expressible by $b_1,b_2$, $G_1$ and $G_2$.  By
\eqref{5.12}  and \eqref{5.14}, every $e\in K_{b_1b_1}$ has the expression
\begin{equation}
e(x,y) = _i\sum ^n_{i,j=1} c_{ij}\xi _i(x)\eta _j(y) + \sum ^n_{i=1} b_2(y)h
''_i(y)\xi _i(x)
+ \sum ^n_{j=1}b_1(x)h'_j(x)_j(y),
\label{eq5.17}
\end{equation}
where, for $i,j = 1,\dots, n$, we have $c_{ij} \in \C $ and $h'_j,
h''_j \in H^2(\T)$ are one-variable functions.  From all the above and
the fact that $e\in K_{b_1b_2}$, we have
\begin{equation}
\thickmuskip=.5 \thickmuskip
\medmuskip=.5 \medmuskip
\!\eqalign{
\Gamma _\phi e =\overline{b_1(x)} \overline{b_2(x)}\biggl(
\sum ^n_{i,j=1}\bigl(G_1(x)G_2(w_j) + G_1(z_i)G_2(y) - G_1(z_i)G_2(w_j)\bigr)
\cdot  \xi _i(x)\eta _j(y)c_{ij} \qquad&\cr+ \sum ^n_{i=1}
G_1(z_i)G_2(y)\xi _i(x)b_2(y)h''_i(y) +
\sum ^n_{j=1}G_1(x)G_2(w_j)b_1(x)\eta _j(y)h'_j(x)\biggr).&}\!
\label{eq5.18}
\end{equation}

\begin{remark}  
Expressions \eqref{5.17}--\eqref{5.18} allow us to check that the
three equivalent conditions of Proposition~\ref{5.2} are also
equivalent to the positive-definiteness of a finite Pick matrix,
defined in terms of $b_1,b_2$, $G_1$ and $G_2$, whose elements are
bounded operators acting in $H^2(\T)$, $L^2(\T)$ or from $H^2(\T)$ to
$L^2(\T)$.
\end{remark}  

Expression \eqref{5.17} shows that $K_{b_1b_2} = K^0 \oplus K^1 \oplus
K^2$, where $K^0$ is the direct sum of the $n^2$ one-dimensional
spaces $\C \xi _i(x)\eta _j(y)$, $K^1$ is the direct sum of the $n$
subspaces $b_2(y)\xi _i(x)H^2(\T)$,
and $K^2$ is the direct sum of the
$n$ subspaces $b_1(x)\eta _j(y)H^2(\T)$.

Whenever $F^1,F^2,F^3,F^4 \in H^\infty (\T)$, it is clear that, for $i,j
= 1,\dots, n$, we have 
$$
\eqalign{(F^1(x) + F^2(y)) \xi _i(x)\eta _j(y) &\in K_{b_1b_2},\cr
F^3(y)\xi _i(x)b_2(y)H^2(\T) &\subset  K_{b_1b_2},\cr
F^4(x)\eta _j(y) b_1(x)H^2(\T) &\subset  K_{b_1b_2}.}
$$

Accordingly, we say that an operator $T: K_{b_1b_2} \rightarrow  K_{b_1b_2}$ is
a multiplier in $K_{b_1b_2}$ if, for $i,j = 1,\dots, n$,
\begin{eqnarray}
\label{eq5.19}
T\xi _i(x)\eta _j(y) &=& (F^1_{ij}(x) + F^2_{ij}(y)) \xi _i(x)\eta _j(y),\\
\label{eq5.19a}
T\xi _i(x)b_2(y)h''(y) &=& F^3_i(y)\xi _i(x)b_2(y) h''(y)\quad\hbox{for $h''
\in H^\infty (\T)$,}\\
\label{eq5.19b}
T\eta _j(y)b_1(x)h'(x) &=& F^4_j(x) \xi _1(x)\eta _j(y) h'(x)\quad\hbox{for $h'
\in H^\infty (\T).$}
\end{eqnarray}

The development above implies the following result:

\begin{theorem}
\label{5.3}
Given two one-dimensional Blaschke products
$b_1$ and $b_2$, with simple zeros at $z_1,\dots,z_n$ and
$w_1,\dots,w_n$, respectively, and given $G_1,G_2
\in H^\infty (\T)$, let $\phi  = (\bar b_1 \otimes  \bar b_2)
(G_1 \otimes  G_2)$. If $\Gamma _\phi $ is the Hankel operator defined
by symbol $\phi $, then $\|\Gamma ^\phi \| = \|\Gamma _\phi \|$,
where $\Gamma ^\phi $ is the multiplier in $K_{b_1b_2}$ $($in the
sense of \eqref{5.19}--\eqref{5.19b}$)$ defined by
$$
\eqalign{
F^1_{ij}(x) + F^2_{ij}(y) &= G_1(x)G_2(w_j) + G_1(z_i) G_2(y) -
G_1(z_i)G_2(w_j),\cr
F^3_i(y) &= G_1(z_i)G_2(y)\cr
F^4_j(x) &= G_1(x) G_2(w_j)}
$$
for  $i,j = 1,\dots,n$.
\end{theorem}

Theorem~\ref{5.3} has a valid formulation in $\T^d$, for $d \geq 1$.  For
$d = 1$ it reduces to the Pick formula.

Theorem~\ref{5.3} allows us to write the boundedness condition
$\|\Gamma_\phi\|\le 1$ as a formula of Pick matrix type, but more
complicated than in the one-dimensional case, and we will not go into
the details here.  Still, remark that the verification of boundedness
of the norm of $\Gamma^\phi$ is not as involved as that for the
restriction of $\Gamma_\phi$ to the model subspace $K_{b_1b_2}$
(condition (b) of Proposition~\ref{5.2}), since it is done through the
defining properties \eqref{5.19}--\eqref{5.19b} of multipliers.


\begin{thebibliography}
\frenchspacing

\bibitem[AAK]{AAK} 
V. M.  Adamjan, V. Z.  Arov and M. G.  Krein, Analytic properties of
Schmidt pairs of a Hankel operator and generalized Schur--Takagi
problem, Mat.  Sbornik {\bf 86} (1971), 33--73. 

\bibitem[ACS]{ACS} 
R.  Arocena, M.  Cotlar and C.  Sadosky, Weighted inequalities in $L^2$ and
lifting   properties, Adv.  Math.  Suppl.  Stud.  {\bf 7A} (1981), 95--128. 

\bibitem[Ag]{Ag}  Jim Agler, Interpolation, J.  Funct.  Anal., to appear. 

\bibitem[Am]{Am}  E. Amar, Les th\'eor\`emes de Schwarz--Pick et
Nevanlinna en plusieurs variables complexes, dans Th\`ese,
Universit\'e de Paris-Sud, Orsay, 1977.

\bibitem[AhC]{AhC} 
P.  Ahern and D. N.  Clark, Invariant subspaces and analytic continuation in
several   variables, J.  Math.  Mech.  {\bf 19} (1969/70), 963--969. 

\bibitem[BH]{BH}
J.  Ball and W.  Helton, A Beurling--Lax theorem for the Lie group
$U(m,n)$ which contains most classical interpolation theory, J.
Operator Theory {\bf 9} (1983), 107--142.

\bibitem[ChF1]{ChF1} 
S.-Y.  Alice Chang and R.  Fefferman, A continuous version of the duality
of $H^1$ and
\BMO/ on the bi-disc, Ann.  of Math.  {\bf 112} (1980), 179--201. 

\bibitem[ChF2]{ChF2} 
S.-Y.  Alice Chang and Robert Fefferman, Some recent developments in
Fourier analysis and $H^p$ theory on product spaces, Bull.  Amer.
Math.  Soc.  {\bf 12} (1985), 1--43.

\bibitem[CLW]{CLW} 
Brian Cole, Keith Lewis and John Wermer, Pick conditions on a uniform
algebra
and von Neumann inequality, J.  Funct.  Anal.  {\bf 107} (1992), 235--254. 

\bibitem[CS1]{CS1} 
M.  Cotlar and C.  Sadosky, The Helson--Szeg\H{o} theorem in $L^p$ of
the bidimensional   torus, Contemp.  Math.  {\bf 107} (1990), 19--37. 

\bibitem[CS2]{CS2}
M.  Cotlar and C.  Sadosky,
Abstract, weighted and multi-dimensional AAK theorems, and the
singular numbers of Sarason commutants, Int.  Eqs.  and Op.  Th.  {\bf
17} (1993), 169--201.

\bibitem[CS3]{CS3}
M.  Cotlar and C.  Sadosky, Nehari and Nevanlinna--Pick problems and
homomorphic extensions in the polydisk in terms of restricted \BMO/,
J.  Funct.  Anal.  {\bf 121} (1994), 205--210.

\bibitem[CS4]{CS4}
M.  Cotlar and C.  Sadosky,
A polydisk version of Beurling's characterization for invariant
subspaces of finite cotype, Proc.  Intl.  Congress on
Complex and Hypercomplex Analysis, Mexico 1994, to appear. 

\bibitem[CW]{CW}
B. Cole and J. Wermer, ``Pick interpolation, von Neumann inequalities
and hyperconvex sets,'' pp.~98--129 in Complex Potential Theory, P. M.
Gauthier (ed.), NATO Adv. Sci. Inst. Ser. C {\bf439}, Kluwer,
Dordrecht, 1994.

\bibitem[HS]{HS}  
H.  Helson and G.  Szeg\H{o}, A problem in prediction theory, Ann. 
Math.  Pure Appl.  {\bf 51} (1960), 107--138. 

\bibitem[N]{N}
Z.  Nehari, On bilinear forms, Ann.  of Math.  {\bf 68} (1957), 153--162. 

\bibitem[Ni]{Ni}
N. K.  Nikolskii, Treatise on the Shift Operator, Springer,
Berlin, 1986. 

\bibitem[P]{P}  
G.  Pick, \"Uber die Beschr\"ankungen analytischer Funktionen, welche
durch vorgegebene Funktionswerte bewirht werden, Math.  Ann.  {\bf 77}
(1961), 7--23.

\bibitem[T1]{T1}  
S.  Treil, The theorem of Adamjan--Arov--Krein: Vector variant, Publ.
Seminar LOMI Leningrad {\bf 141} (1985), 56--72 (in Russian).

\bibitem[T2]{T2}  
S.  Treil, Hankel operators, imbedding theorems and bases of covariant
subspaces of the shift of higher multiplicity, Algebra and Analysis
{\bf 1} (1989), 200--234 (in Russian).
\end{thebibliography}
\end{document}